\documentclass[12pt]{amsart}
\usepackage{amscd}
\newtheorem{thm}{Theorem}
\newtheorem{prop}[thm]{Proposition}
\newtheorem{lem}[thm]{Lemma}
\newtheorem{cor}[thm]{Corollary}

\theoremstyle{remark}
\newtheorem{rem}[thm]{Remark}
\theoremstyle{definition}
\newtheorem{dfn}[thm]{Definition}
\newtheorem{ex}[thm]{Example}

\newcommand{\de}{\delta}

\newcommand{\Hom}{\operatorname{Hom}}

\newcommand{\LL}{\mathcal{ L}}

\newcommand{\N}{\mathbb{ N}}
\newcommand{\Q}{\mathbb{ Q}}
\newcommand{\R}{\mathbb{ R}}

\newcommand{\dg}{\operatorname{deg}} 
\newcommand{\rk}{\operatorname{rk}}
\newcommand{\even}{\operatorname{even}}
\newcommand{\odd}{\operatorname{odd}} 

\newcommand{\Ker}{\operatorname{Ker}}

\DeclareMathOperator{\dummygg}{\mathfrak{g}}
\renewcommand{\gg}{\dummygg}
\DeclareMathOperator{\hh}{\mathfrak{h}}
\DeclareMathOperator{\TT}{\mathfrak{t}}
\DeclareMathOperator{\s}{\mathfrak{s}}

\title{Rational homotopy groups of generalised symmetric spaces}
\author{S.~Terzi\'c}
\address{Mathematisches Institut, Ludwig-Maximilians-Universit\"at 
M\"unchen, Theresienstr.~39, 80333 M\"unchen, Germany}
\email{terzic@mathematik.uni-muenchen.de}
\email{sterzic@cg.ac.yu}
\thanks{The author is 
supported by the {\it DFG Graduiertenkolleg ``Mathematik im Bereich ihrer
Wechselwirkung mit der Physik''} and is a member of EDGE, Research Training Network
HPRN-CT-2000-00101, supported by The European Human Potential Programme.}

\date{\today; MSC2000: 55P62, 22Exx, 22F30}

\begin{document}

\maketitle

\begin{abstract}
We obtain explicit formulas for the rational homotopy groups of generalised symmetric spaces, i.e., the  homogeneous spaces for which the isotropy subgroup appears as the fixed point group of some finite order automorphism of the group.
In particular, this gives explicit formulas for the rational homotopy groups
of all classical compact symmetric spaces.
\end{abstract}

\section{Introduction}

In this paper we consider the problem of calculating the rational homotopy groups of homogeneous spaces $G/H$, where $G$ is a compact connected semisimple Lie group, and 
$H$ is a closed connected subgroup. As a main result, we give explicit formulas for the rational homotopy groups in the case when $H$ is the fixed point subgroup of some finite order automorphism of $G$. We call such  homogeneous spaces generalised symmetric spaces, and, obviously, they contain the symmetric spaces. In some places in the literature these spaces are known as $k$ - symmetric spaces. 

The problem of computing rational homotopy groups was theoretically solved in the 1970's
by Sullivan's minimal model theory~\cite{S}. One of the results of this theory is that to each connected commutative graded differential algebra  $(A, d_{A})$ one can assign minimal model (uniquely up to isomorphism), i.e., a free commutative graded differential algebra $(\mu _{A}, d)$  with decomposable coboundaries which is quasi-isomorphic to the algebra $(A, d_{A})$. In particular to each smooth manifold $M$, one can assign a minimal model. In the case when $M$ is a homotopy simple manifold of a finite type its minimal model completely classifies its rational homotopy type. Under the same assumptions on $M$, the minimal model  contains complete information on the rational homotopy groups of  the manifold~\cite{BG},~\cite{lehm}. More precisely, the degrees and multiplicities of its  generators determine the non-trivial rational homotopy groups of $M$.

Nevertheless, the computation of rational homotopy groups remains quite a difficult problem in general.

There are still very few homogeneous spaces for which all
rational homotopy groups are known. Among them are those spaces with free cohomology algebras to which the Cartan-Serre theorem can be applied directly. Examples of these are the semisimple Lie groups and the spheres of odd dimensions. 
The cohomology rings of the spheres of even dimension and of the complex projective spaces are not free, but they have  quite simple real cohomology algebras and the application of Sullivan's minimal model theory immediately gives their rational homotopy groups~\cite{S}.

Further examples are the homogeneous spaces $G/S$, where $S$ is a torus in a semisimple Lie group $G$. Their rational homotopy groups can be immediately obtained by applying the above results to the exact homotopy sequence of
the corresponding fibration.

The main difficulties in the application of Sullivan's minimal model theory are in connection with the explicit construction of the minimal model of certain connected commutative graded differential algebras.
 
This problem  is solved in the case when  
$(A, d_{A})$ is a free simply connected commutative 
graded  differential algebra. Namely, for such an algebra, Sullivan's minimal model theory gives a method which says that the generators of its minimal model (which in general are not canonical) correspond to those generators of $A$ whose differential has a non-linear part, and which do not appear as a linear part in the coboundaries. In order to apply this 
method, one needs to have explicit expressions for the action of the differential $d_{A}$ on the generators  of $A$. 

Note that Sullivan's minimal model theory, as well as many other results related to rational homotopy theory, are explained in detail in the recent book~\cite{FHT}.

Let $G$ be a compact connected Lie group and $H$ a connected closed subgroup.
It is a well known important fact that to the homogeneous space $G/H$ one can  assign a Cartan algebra~\cite{Borel} which  turns out to be weakly equivalent to its de Rham algebra of differential forms~\cite{lehm}. Thus, the minimal model of $G/H$ is the minimal model of its Cartan algebra. If we assume $G$ to be semisimple, the Cartan algebra of $G/H$ becomes simply connected and,
in general, one can apply Sullivan's method for the computation of its minimal model. In practice, to apply this method we need to know
how  the differential in the Cartan algebra acts, or, what is equivalent, the inclusion of the maximal Abelian subalgebra for the subgroup $H$ into the maximal Abelian subalgebra of the group $G$.

In this paper, from the knowledge of these inclusions in the case of 
generalised symmetric spaces~\cite{Helgason},~\cite{Kac1},~\cite{T},~\cite{T1} we derive the rational homotopy
groups of these spaces.

The paper is organised as follows. In Section~\ref{MM} we review the basic definitions and theorems from minimal model theory that we will need. Section~\ref{HS} contains results on the rational homotopy groups of homogeneous spaces of semisimple compact Lie groups $G$, which we can get by direct applications of Sullivan's algorithm. In Section~\ref{GSS} we state and prove the main theorem on rational homotopy groups of
generalised symmetric spaces. Then, in Appendix A, by applying the main theorem, we give   
the table for the  rational homotopy groups of irreducible simply connected compact symmetric 
spaces. 

{\sl Acknowledgment}: I would like to thank Dieter Kotschick for the support.  

\section{ On minimal model theory}\label{MM} 
Let $(A, d_{A})$ be a commutative graded differential algebra over a filed $k$ of characteristic zero. By definition, a differential graded algebra $(\mu _{A}, d)$ is called a model for $(A, d_{A})$ if there exists a differential graded algebras morphism
\[
h_{A} : (\mu _{A}, d)\to (A, d_{A})
\]
inducing an isomorphism in cohomology. Such a morphism is called quasi-isomorphism and the differential graded algebras for which it exists are called quasi isomorphic or weakly equivalent.

 If $(\mu _{A}, d)$ is a free algebra in a sense that $\mu _{A} = \wedge V$, for a graded vector space $V$, then it is called a free model for $(A, d_{A})$. The notation $\mu _{A} = \wedge V$ means 
that, as a graded vector space $\mu _{A}$ is polynomial algebra on even degree elements $V^{even}$ and an exterior algebra on odd degree elements $V^{odd}$.

\begin{dfn}
A differential graded algebra $(\mu _{A}, d)$ is called a minimal model of $(A, d_{A})$ if:
\begin{itemize}
\item $(\mu _{A},d) = (\wedge V, d)$ is a free model for $(A, d_{A})$ \ ;
\item $d$ is indecomposable in the following sense: for a fixed set $V =\{ P_{\alpha}, \alpha \in I\}$ of free generators of $\mu _{A}$, where $I$ is a poset, we have that for any  $P_{\alpha}\in V$, $d(P_{\alpha})$ is a polynomial in generators $P_{\beta}$, $\beta <\alpha$, with no linear part.
\end{itemize}
\end{dfn}
Then, it can be proved (~\cite{lehm},~\cite{GM},~\cite{FHT}) that for any connected commutative differential $\N$-graded algebra $(A, d_{A})$ ($H^{0}(A, d_{A}) = k$) there exists minimal model (unique up to isomorphism). 

By definition minimal model of some  smooth (connected) manifold $M$ will be minimal model of its  de Rham algebra $\Omega ^{*}(M)$ of differential forms. 
Since $\Omega ^{*}(M)$ is connected, we have the existence of minimal model, $\mu (M)$, for any such a manifold. 
It turns out that, under suitable conditions, the knowing of the minimal model of a manifold permits us to compute the ranks of its homotopy groups. 

We say that the manifold $M$ is a homotopy simple if $\pi _{1}(M)$ is commutative and acts trivially on $\pi _{p}(M)$, for $p\geq 2$.
Denote by $\mu ^{++}(M)$ the set of all indecomposable elements in $\mu (M)$, i.e., 
$\mu ^{++}(M) = \mu (M)/\mu ^{+}(M)\cdot \mu ^{+}(M)$.

Then, in~\cite{BG},~\cite{S} the following theorem is proved.

\begin{thm} \label{thm:ranks} Let M be a homotopy simple manifold of
finite type. Then $\mu ^{++} (M)_{r}\cong \Hom (\pi _{r}(M), {\R})$, $r\geq 1$. 
Thus, $\rk \pi _{r}(M) = \dim \mu ^{++} (M)_{r}$, $r\geq 1$.
\end{thm}
\begin{rem}\label{st}
Any homogeneous space of a connected Lie group, whose stabilizer of a point is
connected, is homotopy simple (~\cite{ST}, $\S 16$).
\end{rem}
\begin{rem}
The above theorem is true for all nilpotent spaces of a finite type, see~\cite{BG},~\cite{lehm}.
By Remark~\ref{st}, Theorem~\ref{thm:ranks} is sufficient for our purposes.
\end{rem} 

To apply the above theorem to compute the rational homotopy groups of some manifold one needs an explicit expression for the algebraic part  of its minimal model ( obviously, we do not need information on its differential). In general it is not clear how the algebraic part of the minimal model for certain algebra or  manifold can be computed effectively. We will use the algorithm provided by Sullivan which allows us to compute the minimal model of any 
free commutative graded differential algebra generated by elements of degree $>1$. In particular, we will see later that, under suitable conditions, this algorithm yields the minimal model of certain compact homogeneous spaces.

\subsection{Minimal model of a simply connected free CGDA}~\label{alg}

We say that  the commutative graded differential algebra $(A, d_{A})$ is simply connected if $H^{0}(A, d_{A}) = k$ and $H^{1}(A, d_{A}) = 0$. 
We will describe the algorithm given by Sullivan (~\cite{lehm},~\cite{FHT}) for the computation of the minimal model of a simply connected free graded differential algebra.

Let $(\wedge V, d)$ be a free graded  differential algebra. Assume that it is simply connected. Then
$V = \oplus _{k\geq 1} V^{k}$ and $H^{1}(\wedge V, d) = 0$. Denote by $\wedge ^{++}V$ the ideal in $\wedge V$
generated by all decomposable elements. Let $\de : V^{n}\to V^{n+1}$ be a
differential given by the formula
\[
\de = \pi\circ d \ ,
\]
where $\pi$ is the canonical projection
\[
\pi : V^{n+1}\oplus (\wedge ^{++}V)^{n+1}\to V^{n+1} \ .
\]
Let $V^{'}$ be a complement of $\Im \de$ in $\Ker \de$ and denote by $W$
a complement of $\Ker \de$ in $V$. Thus, 
\[
V = \Im \de\oplus V^{'}\oplus W \ .
\]
By straightforward calculations one can  obtain
\[
\wedge V = \wedge (V^{'})\otimes \wedge (W^{'}\oplus W) \ ,
\]
where $W^{'} = d(W)\subset \wedge V$. Let $C = \wedge (W^{'}\oplus W)$.
Denote by $\langle C^{+}\rangle $ the ideal in $\wedge V$ generated by the elements
of degree $>0$ in $C$. Then the canonical projection
\[
\rho : (\wedge V, d)\to (\wedge V/\langle C^{+}\rangle , d_{1})
\]
is a quasi-isomorphism of differential graded algebras ($d_{1}$ is a  
differential induced by $d$). Moreover, $\wedge V/\langle C^{+}\rangle \cong
\wedge V^{'}$ and, therefore, is free. Being simply connected, it is also minimal.

\section{ On rational homotopy groups of compact homogeneous spaces}\label{HS}

\subsection{General theory}\label{GT}
Let $G$ be a compact connected Lie group and $H$ a closed connected
subgroup such that  $\rk G = n$ and $\rk H = r$. By $\TT$ and $\s$ we denote the maximal Abelian subalgebras for $G$ and $H$ respectively.
Let us consider the homogeneous space $G/H$.


The real Cartan algebra $(C, d)$ of  $G/H$ has the form $C = H^{*}(B_H)\otimes H^{*}(G)$. We recall some facts from~\cite{Borel} concerning the structures of $H^{*}(B_H)$ and $H^{*}(G)$.  
The cohomology algebra of the classifying 
space $B_H$ is
the algebra of the Weil invariant polynomials on
$\s$. This algebra allows $r$ generators $Q_{l_1},\ldots ,Q_{l_r}$ whose degrees are $\dg Q_i = 2l_i$, $1\leq i\leq r$, i.e., $H^{*}(B_{H})\cong \R [Q_{l_1},\ldots ,Q_{l_r}]$.
Also $H^{*}(G)$ is the exterior algebra over $n$ universal transgressive elements $z_{k_1},\ldots ,z_{k_n}$ whose degrees are $\dg z_{k_i} = 2k_i - 1$, $1\leq i\leq n$, i.e., $H^{*}(G) = \wedge (z_{k_1},\ldots ,z_{k_n})$. By $P_{k_1},\ldots ,P_{k_n}$, $\dg P_{k_i}=2k_{i}$, we denote the generators for the cohomology algebra $H^{*}(B_G)$ of the classifying space for $G$ corresponding to $z_{k_1},\ldots , z_{k_n}$ by transgression in the universal bundle for $G$. The numbers $k_1,\ldots ,
k_n$ ($l_1, \ldots ,l_r$) are often called the exponents of a group $G$ ( subgroup $H$). For the simple compact Lie groups these numbers, as well as the generators of the Weil invariant polynomial algebras, are well known ~\cite{MT}. We give their list at the Appendix~\ref{B}, since we will need their explicit expression. Further in the paper, by $\nu (k_{i})$ ($\nu (l_{j})$) we will denote the multiplicity of the exponent $k_{i}$ ($l_{j}$).

 Differential in the Cartan algebra is defined as 
$d(b\otimes 1) = 0$, $b\in H^{*}(B_H)$ and 
$d(1\otimes z_i) = \rho ^{*}(P_i)\otimes 1$. We denoted by $\rho^{*}$ the restriction of the Weil invariant polynomial algebra  $\mathbb{ R}[\TT ]^{W_G}$ to the Weil invariant polynomial algebra $\mathbb{ R}[\s ]^{W_H}$, see~\cite{Borel}.

Cartan theorem~\cite{Borel} states that the real cohomology algebra of
the homogeneous space $G/H$  is isomorphic to the cohomology algebra of its Cartan
algebra. But, it turns out that more is true~\cite{Onishchik}:

\begin{thm}
The Cartan algebra of the homogeneous space $G/H$ is weakly equivalent 
to its de Rham algebra of differential forms.
\end{thm}

Thus, the minimal model of the homogeneous space $G/H$ is the minimal model
of its Cartan algebra. 
\begin{rem}
The above theorem implies well known fact that the simply connected compact homogeneous spaces are rationally
elliptic~\cite{FHT}.
\end{rem}  
Since the Cartan algebra is a free
graded differential algebra, under assumption on a group $G$ to be semisimple, we can apply algorithm given in~\ref{alg} for 
computing its minimal model.
\begin{rem}
Note that in the case when $\rk G=\rk H$, we can always assume $G$ to be semisimple, according to the results on decomposition of $G$ and $H$, see~\cite{Onishchik}(p. 61).
\end{rem}
In the notation of section~\ref{alg} we have 
\begin{equation} \label{primjena}
V = \LL (Q_1,\ldots ,Q_r)\oplus \LL (z_1,\ldots ,z_n) \ ,
\end{equation}
and differential $\de : V^{n}\to V^{n+1}$ acts as follows:
\begin{itemize}
\item $\de (Q_i) = 0$, $1\leq i\leq r$
\item $\de (z_i)$ is a linear part of $\rho ^{*}(P_i)$ in
$\mathbb{ R}[\s ]^{W_H}$, $1\leq i\leq n$.
\end{itemize}
Let $\hh$ be a Lie algebra for $H$ . Then we have obvious decomposition
\[
\hh = \xi (\hh )\oplus \hh ^{'} \ ,
\]
on its center $\xi (\hh )$ and its semisimple part $\hh ^{'}$. This implies
\[
\s = \xi (\hh)\oplus \s ^{'} \ ,
\]
where $\s ^{'}$ is the maximal Abelian subalgebra in $\hh ^{'}$. 

Thus, 
\[
\mathbb{ R}[\s ]^{W_{\hh }} = \mathbb{ R}[\xi (\hh )]\otimes \mathbb{ R}[\s ^{'}]^{W_{\hh ^{'}}} \ .
\]
It follows that the above differential $\de$ can be written in the form 
$\de = \de _{1} + \de _{2}$, where $\de _{1}$ and $\de _{2}$ are the linear parts of $d$ in $\R [\xi (\hh )]$ and $\R [\s ^{'}]^{W_{\hh ^{'}}}$ respectively.
 
Now, $\R [\xi (\hh )]$ is a polynomial algebra of $\dim \xi (\hh )$ variables whose all degrees are $2$, and since $G$ is a semisimple, all the generators in $H^{*}(B_{G})$ are of degree $\geq 4$. It follows that 
$\de _{1} = 0$,
i.e., $\de = \de _{2}$. 

Thus, the computation of the  algebraic part ($\mu = V^{'}$) of the minimal model for $G/H$ is equivalent to the computation of the kernel and  the image of $\de$ restricted to $\s ^{'}$, i.e., we need to consider
 \begin{equation}\label{restr}
\de : \LL (z_1,\ldots ,z_n) \to \LL (Q_{i_1},\ldots ,Q_{i_j}) \ ,
\end{equation} 
where $Q_{i_1},\ldots Q_{i_j}$ are the generators in $\R [\s ^{'}]^{W_{\hh ^{'}}}$.

Now, being a semisimple Lie algebra, $\hh ^{'}$ splits into the sum of the simple compact
Lie algebras, i.e.,  $\hh ^{'} = \hh _{1}^{'}\oplus \ldots \oplus \hh  _{q}^{'}$ and, therefore
\begin{equation} \label{ti}
\R [\s ^{'}]^{W_{\hh ^{'}}} = \R [\s _{1}^{'}]^{W_{\hh _{1}^{'}}}\otimes
\ldots \otimes \R [\s _{q}^{'}]^{W_{\hh _{q}^{'}}} \ , 
\end{equation}
where $\s _{1}^{'},\ldots ,\s _{q}^{'}$ are the maximal Abelian subalgebras for 
$\hh _{1}^{'},\ldots , \hh _{q}^{'}$ respectively. It follows, as above, that $\de = \de _{\s _{1}}+\ldots + \de _{\s _{q}}$, where $\de _{\s _{i}}$ is a linear part of the restriction
 \[
\rho ^{*}_{\s _{i}} : \R [\TT ]^{W_{\gg}}\to \R [\s _{i}]^{W_{\hh _{i}}} \ , 
\]
which is induced by $\rho ^{*}$.
Then, obviously, 
$\de (z_{i})\neq 0$ if and only if $\de _{\s _j}(z_{i})\neq 0$ for some $j$, $1\leq j\leq q$. 

\subsection{Some direct applications.}
Let us first illustrate some straightforward applications of the above algorithm  which will help us to determine which are the rational homotopy groups of the homogeneous spaces difficult to compute.  
 
The statements~\ref{flag} to~\ref{thm:odd-even} can be also proved by considering the long exact homotopy sequence of the fibration $H\to G/H\to G$. 
\begin{ex}\label{flag}
Let $G$ be a compact connected semisimple Lie group and $S$ a toral subgroup
of dimension $r$. Then $\R [\s ]^{W_S} = \R [x_1, \ldots ,x_r]$, where $\dg x_i = 2$. Since $G$ is semisimple, $\dg \rho ^{*}(P_i) \geq 4$ for any Weil invariant generator $P_{i}$  in $\R [\TT ]^{W_G}$, so we get that $\de$ vanishes identically. 
Therefore, the algebraic part in the minimal model for $G/S$ coincides with its Cartan algebra:
\[ 
\mu (G/S) = R[x_1,...,x_r]\otimes \wedge (z_1,...,z_n) \ .
\]
Thus, from~\ref{alg} follows:
\[
 \dim \pi _{2}(G/S)\otimes \Q = \dim S, \;\;  \pi _{2p}(G/S)\otimes \Q = 0, \;\; \mbox{for} \;\;
p\geq 2 \ , 
\]
\[
\pi _{2p-1}\otimes \Q = 0, \; \mbox{for} \; p\notin \{ k_1,\ldots ,k_n \} \ ,
\]
\[ 
\dim \pi _{2k_i-1}(G/S)\otimes \Q = \nu (k_i) , \; \mbox{for} \;\; i=\bar{1,n} \ .
\]
\end{ex}
\begin{rem}
Note that, for any homogeneous space $G/H$ such that all $\rho ^{*}(P_{i})$ are decomposable, its minimal model also coincides with its Cartan algebra.
\end{rem}
We can push up the above further. 
\begin{lem}
Let $G$ be a compact connected semisimple Lie group and $H$ a closed connected subgroup. Then 
$\dim \pi _{2}(G/H)\otimes \Q = \dim \xi (H)$.
\end{lem} 
\begin{proof}
Elements of degree $2$ in $\Ker \de$ may appear only if $H$ has non trivial
center $\xi (H)$. Since $\Im \de$ contains no elements of degree two, we get the statement.
\end{proof}
Similarly, using~\eqref{primjena}, one can easily prove the following. 
\begin{prop} \label{thm:odd-even}
Let $G$ be a semisimple compact connected Lie group, $H$ a connected closed subgroup, 
$k_1,\ldots ,k_n$ the  exponents of the group $G$ and $l_1,\ldots ,l_r$ the
exponents of the subgroup $H$. Then
\begin{enumerate}
\item $ \pi _{2p}(G/H)\otimes \Q = 0$ for $p\notin \{l_1,\ldots ,l_r\}$;
\item $\pi _{2p-1}(G/H)\otimes \Q = 0$ for $p\notin \{k_1,\ldots ,k_n\}$;
\item $\dim \pi _{2k_i - 1}(G/H)\otimes
\Q  = \nu (k_i)$ for  $k_i\notin \{l_1,\ldots ,l_r\}$ , $1\leq i\leq n$.
\item  $\dim \pi _{2l_i}(G/H)\otimes
\Q = \nu (l_i)$ for
$l_i\notin \{k_1,\ldots ,k_n\}$, $1\leq i\leq r$.
\end{enumerate}
\end{prop}
 
From the above proposition immediately follows:
\begin{cor}
$\dim \pi _{\even}(G/H)\otimes \Q \leq \rk H$
and $\dim \pi _{\odd}(G/H)\otimes \Q \leq \rk G$.
\end{cor}
\subsection{Further simple applications.} Concerning Proposition~\ref{thm:odd-even}, in order to compute the rational homotopy groups of the 
homogeneous space $G/H$, it remains to see what are the indecomposable elements 
in $d(z_{k_i})$, for $z_{k_i}\in H^{*}(G)$ and  $k_{i}\in \{l_1,\ldots , \l_r\}$.
 
Note that, when $\nu (k_{i}) = 1$ 
we do not need to know how the indecomposable part of $d(z_{k_i})$ looks like. When this is the case, we are only interested in the vanishing of $\de (z_{k_i})$. More precisely, the following statement is obvious.
\begin{prop}\label{nu=1}
Let $k_{i} = l_{j}$ be the common exponent for a semisimple compact connected Lie group $G$ and its closed connected subgroup $H$, such that $\nu (k_{i}) = 1$ and $z_{k_i}\in H^{*}(G)$, $\deg z_{k_i} = 2k_{i}-1$. Then
\begin{itemize}
\item $\dim \pi _{2l_{j}}(G/H)\otimes \Q = \nu (l_{j})$ and $\dim \pi _{2k_{i}-1}\otimes \Q = 1$ if $d^{'}(z_{k_i}) = 0$
\item  $\dim \pi _{2l_{j}}(G/H)\otimes \Q = \nu (l_{j}) - 1$ and $\pi _{2k_{i}-1}\otimes \Q = 0$ if $d^{'}(z_{k_i})\neq 0$.
\end{itemize}
\end{prop}
%

The above is always the case when $G$ is a simple compact Lie group of a type different from $D_{2n}$. If $G$ is of a type $D_{2n}$, then it has exponent $2n$
of multiplicity $2$, while all the others are of multiplicity $1$. This implies the following.
\begin{cor} \label{p:prva}
For any simple compact connected Lie group $G$ that is not of type $D_{2n}$ it follows that  $\rk \pi _{2k-1}(G/H)\leq 1$.
For a group of a type $D_{2n}$, we have $\rk \pi _{4n-1}(G/H)\leq 2$ and $\rk \pi _{2k-1}(G/H)\leq 1$ for $k\neq 2n$.
\end{cor}

\begin{cor}\label{pi34}
Let $G$ be a simple compact Lie group and $H\neq \{e\}$ a closed connected subgroup that is not a torus. Then
\[
\pi _{3}(G/H)\otimes \Q = 0, \; \mbox{and} \;
 \dim \pi _{4}(G/H)\otimes \Q = \nu (l_1=2)-1
\ .
\]
\end{cor}
\begin{proof}
Since $G$ is simple it has only one Weil invariant generator $P_2$ of degree $4$. Concerning the formulas for the Weil invariant generators,~\eqref{agen} -~\eqref{e6gen} and~\cite{MT}, we have that the restriction of $P_{2}$ on any
subspace of $\TT$ is not trivial.
Being of the smallest even degree, it means that $\de _{s_j}(z_2)\neq 0$, for any $j$, $1\leq j\leq q$ defined in~\eqref{ti}. Then, Proposition~\ref{nu=1} 
implies the statement.
\end{proof}
For some restrictive class of homogeneous spaces we can say some more.
\subsection{ Non cohomologous to zero homogeneous spaces.} Recall that, it is said that, the subgroup $H$ is totally non cohomologous to zero
in the group $G$ if the restriction $\rho ^{*}:\R [\TT ]^{W_G}\to \R [\s ]^{W_H}$ is surjective. 

This property easily implies the following.
\begin{lem}
Let the subgroup $H$ be totally non cohomologous to zero in the group $G$. Then
all even rational homotopy groups of $G/H$ are trivial and
\[
\pi _{2k-1}(G/H)\otimes \Q  = \nu (k_{i}) - \nu (l_{j}) \ ,
\]
for any common exponent $k = k_{i} = l_{j}$ for the group $G$ and subgroup $H$. 
\end{lem}
\begin{proof}
Surjectivity of $\rho ^{*}$ implies surjectivity of $\de$. Thus, any exponent of $H$ is also an exponent for $G$, and $\dim (\Im \de )_{2k} = \nu (l_{j})$ for any common exponent $k = k_{i} = l_{j}$. It follows that all even rational homotopy groups are trivial and
$\pi _{2k-1}(G/H)\otimes \Q = \nu (k_{i}) - \nu (l_{j})$.
\end{proof}
\begin{rem}
In this case $H^{*}(G/H)$ is an exterior algebra and the Cartan-Serre theorem can be also applied.
\end{rem}
\begin{rem}
Surjectivity of $\de$ implies surjectivity of $\rho ^{*}$. Thus, the opposite of the previous lemma is also true. If all even rational
homotopy groups of $G/H$ are trivial, then $H$ is totally non cohomologous to zero
in $G$. 
\end{rem}
\subsection{Cartan pair homogeneous spaces.}
Recall~\cite{GHV} that we say that the homogeneous space $G/H$, $\rk G = n$, $\rk H = r$, is  a Cartan pair if we can choose $n$ algebraic
independent generators $P_1,\ldots ,P_n$ of $\mathbb{ R}[\TT ]^{W_G}$ such that 
$\rho ^{*}(P_{r+1}),\ldots ,\rho ^{*}(P_n)$ belong to the ideal in 
$\R [\s ]^{W_H}$ which is generated by $\rho ^{*}(P_1),\ldots ,\rho ^{*}(P_r)$. Furthermore, when this is the case, one can choose $P_{r+1},\ldots ,P_{n}$ such that $\rho ^{*}(P_{r+1})=\ldots =\rho ^{*}(P_{n})=0$, see~\cite{Doan Kuin'}. This immediately implies the following Lemma. 
\begin{lem}
Let $G/H$ be a Cartan pair homogeneous space and let $\dg P_i = 2p_{i}$, $r+1\leq i\leq n$. Then
$\pi _{2p_i - 1}(G/H)\otimes \Q$ are not trivial for all $r+1\leq i\leq n$.
\end{lem}
Obviously, if the subgroup $H$ is totally non homologous to zero in the group $G$, then $G/H$ is a Cartan pair homogeneous space.
\begin{rem}
The previous Lemma is also a direct consequence of the results on cohomology structure and formality of Cartan pair homogeneous spaces.
Is is known~\cite{Onishchik} that the homogeneous space $G/H$ is formal if and only if it is a Cartan pair homogeneous space. For such space it is known~\cite{Onishchik} that its
cohomology algebra has the form
\begin{equation}\label{CC}
H^{*}(G/H)\cong \mathbb{ R}[\s ]^{W_H}/\langle \rho ^{*}(R[\TT ]^{W_G})\rangle \otimes \wedge (z_{r+1},\ldots ,z_n) \ .
\end{equation} 
Thus, the minimal model of $G/H$ is the minimal model of the above cohomology
algebra endowed with the differential equal to zero. 
\end{rem}
\begin{rem}
Note that the spaces having the cohomology algebras of the form~\eqref{CC} generalise the notion of the spaces with good cohomology in terminology of~\cite{BG},~\cite{SMITH}. We mean generalisation in a sense that the construction of the minimal models for the spaces with good cohomology~\cite{BG} can be directly generalised to the spaces with the cohomology algebras of the form~\eqref{CC}.
\end{rem} 
 We see that a lot of the rational homotopy groups of the homogeneous space $G/H$ can be obtained comparing the exponents of $G$ and $H$.
Problems arise for those rational homotopy groups corresponding to the common exponents, since for them we need to 
know differential $d$ in the Cartan algebra for $G/H$. 
Obviously for any such exponent
one can consider  corresponding odd ( taking it as an exponent of the group $G$) and even (taking it as an exponent of the subgroup $H$) rational homotopy group.

In this direction we will obtain results about special type of homogeneous spaces, which
we will call generalised symmetric spaces.

\section{Rational homotopy groups of generalised symmetric spaces}\label{GSS}

\subsection{Generalised symmetric spaces.}
Several generalisations of the symmetric spaces are known and they have been studied
by many authors, see~\cite{Kac},~\cite{Kow},~\cite{Selberg},~\cite{V-F},~\cite{W-G},...
Following~\cite{Kow}, the most general definition for generalising the notion of symmetric spaces by omitting the condition on involutivity is as follows. 
\begin{dfn}
A generalised symmetric space of order $m$ is a triple $(G, H, \Theta)$,
where $G$ is
a connected Lie group, $H\subset G$ is a closed subgroup, and $\Theta$ is
an automorphism of finite order $m$ of the group $G$ satisfying
\[ 
G_{0}^{\Theta}\subseteq H\subseteq G^{\Theta} \ ,
\]
where $G^{\Theta}$ is the fixed point set of $\Theta$ and $G_{0}^{\Theta}$
is its identity component.
\end{dfn}
If we assume $G$ to be semisimple and simply connected, generalised symmetric spaces are given by the triples $(G,G^{\Theta},\Theta )$~\cite{Vinberg},~\cite{OV}.
 Under the same 
assumptions, there is also a bijection between generalised symmetric spaces and
generalised symmetric algebras $(\gg , \gg ^{\Theta}, \Theta )$, where the second ones are defined on the same way 
~\cite{Vinberg},~\cite{OV}.
So, further we will make no difference between  generalised symmetric spaces and general symmetric algebras.

From the classification theory, developed in~\cite{Kac} for the classification of finite order automorphisms of the semisimple Lie algebras
, it follows that, for $\gg$ simple, to any finite order automorphism $\Theta$ we can assign a  sequence $(k, (s_0,\ldots ,s_r), m)$, called the type of automorphism $\Theta$, which completely determines it (up to conjugation by some automorphism 
of $\gg$). Here $k = 1, 2, 3$, then $s_0,\ldots ,s_r$ are 
non-negative integers, without non-trivial common factor, corresponding to the vertices of the Kac-Dynkin diagram $\gg ^{(k)}$ and  $m$ is an order  of $\Theta$. For generalised symmetric algebras this classification implies the following~\cite{Kac1},~\cite{Helgason}.
\begin{enumerate} \label{clas}
\item Any generalised symmetric algebra $(\gg, \gg ^{\Theta}, \Theta)$ can be decomposed as 
\begin{equation}\label{decomp}
\gg ^{\Theta} = \TT ^{r-t}\oplus \gg ^{'} \ ,
\end{equation}
where $\gg ^{'}$ is a semisimple Lie algebra whose Dynkin diagram is subdiagram 
of the Kac-Dynkin diagram $\gg ^{(k)}$ consisting of the vertices $\{ i_1,\ldots ,i_t\}$ given by the conditions $s_{i_1} = \ldots = s_{i_t} = 0$.
\item
Generalised symmetric algebras $(\gg , \gg ^{\Theta}, \Theta )$ and $(\gg, \gg ^{\tilde {\Theta}}, \tilde {\Theta})$ are isomorphic if and only if $k = \tilde {k}$ and
the sequence $(s_0,\ldots ,s_r)$ can be transformed into the sequence $(\tilde {s}_{0},\ldots ,\tilde {s}_{r})$ by some automorphism of the Kac-Dynkin diagram
$\gg ^{(k)}$.
\end{enumerate} 
In this way, using above results, an explicit list of all generalised symmetric spaces of compact simple simply connected Lie groups is given in~\cite{JV}.  

\subsection{The main result.} 
Our main result is the following theorem giving explicit formulas for the rational homotopy groups of generalised symmetric spaces. 
\begin{thm}\label{main}
Let $G/H$ be  a generalised symmetric space, where $G$ is a simple 
compact Lie group. Let further $\gg$ be the Lie algebra for $G$, $\gg ^{'}$ the semisimple part in the Lie algebra for $H$, $\{ k_1,\ldots ,k_n \}$ the exponents 
of the group $G$, $\{ l_1,\ldots ,l_r \}$  the exponents of the subgroup $H$ and
$p = k_{i} = l_{j}$.
Then:
\begin{enumerate}
\item  If $\rk H = \rk G$ we have 
\begin{itemize}
\item 
$\dim \pi _{2p}(G/H)\otimes \Q = \nu (l_{j}) -1$\ ,
\item
$\pi _{2p-1}(G/H) \otimes \Q = 0$ for $(\gg, p)\neq (D_{2n}, 2n)$\ ,
\item 
$\pi _{4n-1}(G/H)\otimes \Q = \Q$ for $\gg = D_{2n}$ and $\gg ^{'}\neq A_{2n-1}$\ ,\\
$\pi _{4n-1}(G/H)\otimes \Q = 0$ for $(\gg, \gg ^{'}) = (D_{2n}, A_{2n-1})$\ .
\end{itemize}
\item If $\rk H <\rk G$ we have
\begin{itemize}
\item 
for $p$  odd :\\
$\dim \pi _{2p}(G/H)\otimes \Q = \nu (l_{j})$,  $\pi_{2p-1}(G/H)\otimes \Q = \Q$\ ,
\item
for $p$ even: \\
 $\dim \pi_{2p}(G/H)\otimes \Q = \nu (l_{j}) -1$\ , \\
$\pi _{2p-1}(G/H)\otimes \Q = 0$ for  $(\gg ,p)\neq (D_{2k}, 2k)$\ ,\\
$\pi _{2p-1}(G/H)\otimes \Q = \Q$ for $(\gg ,p) = (D_{2k}, 2k)$\ .
\end{itemize}
\end{enumerate}  
\end{thm}
\begin{rem}
In~\cite{T} it was proved that all generalised symmetric spaces of  simple compact Lie groups are  Cartan pair homogeneous spaces. Therefore~\cite{KT} they are formal in the sense of Sullivan.
\end{rem}
Before going to prove the above theorem, let us note one example showing that the theorem can not be extended to the Cartan pair homogeneous spaces. 
\begin{ex}
All compact Riemannian homogeneous spaces with an irreducible isotropy group are classified by O.~V.~Manturov~\cite{OVM}. 
Let us consider such a space $G/H$ given by
$\gg = A_{N-1}$ and $\hh = A_{n-1}$, where $N = \frac{n(n-1)}{2}$, $n\geq 6$. In~\cite{ Doan Kuin'} is shown that the restriction of the Weil invariant generators $P_{a} = x_1^{a}+\ldots +x_{N}^{a}$, $a = 2,3,\ldots ,N$ of $A_{N-1}$, to the maximal Abelian 
subalgebra of $A_{n-1}$ is given by
\[
\tilde {P_a} = \sum _{i<j}(y_i + y_j)^{a} = (n-2^{a-1})p_a + \frac{1}{2}\sum_{s=2}^{a-2}c_{a}^{2}p_{a-s}p_{s} \ ,
\]
where $p_s = y_1^s +\ldots +y_n^s$, $s = 2,\ldots ,N$.

It is also proved that for $n = 2^{q-1}$, instead of standard $Q_2,\ldots ,Q_n$(where $Q_{i}=p_{i}$ for $2\leq i\leq n$) we can take $\tilde {P}_{2}\ldots , \tilde {P}_{q-1}, Q_q, \tilde {P}_{q+1},\ldots ,\tilde {P}_{n}$ for the generators in $\R [\s ]^{W_H}$. When this is the case, the restriction $\rho ^{*}$ in the Cartan algebra for $G/H$ is given by:
$$
\rho ^{*}(P_{a})=\tilde {P_a}, \; 2\leq a\leq n, a\neq q, \;\; 
\rho ^{*}(P_{q})=0 \ ,
$$
$$ 
\rho ^{*}(P_{kq})=\alpha _{k}Q_{q}^{k}, \; n< kq\leq N, \;\; \rho ^{*}(P_{a})=0, \; n< a \leq N, a\neq kq \ .
$$
Since $k>1$, it follows that the differential  $\de$ is given by 
$$
\de(P_{a})=\tilde {P_a}, \; 2\leq a\leq n, a\neq q \ ,
$$
$$ 
\de (P_{q})=0 \;\; \mbox{and} \;\; \de (P_{a})=0, \; n< a \leq N, \ .
$$

Thus, $V^{'}=\LL (Q_{q}, z_{q}, z_{n+1},\ldots ,z_{N})$ and since $\dg Q_{q}=2q$ and $\dg z_{i}=2i-1$, we obtain that
\[
\rk \pi _{2q}(G/H)=\rk \pi _{2q-1}(G/H)=1 \ ,
\]
\[
\rk \pi _{2i-1}(G/H)=1, \; n+1\leq i\leq N \ ,
\]
\[
\rk \pi _{p}(G/H)=0, \; p\neq 2q, 2q-1, 2i-1, \; n+1\leq i\leq N \ .
\]
 On the other hand, for any odd  common exponent $p$ for $A_{N-1}$ and $A_{n-1}$, (i.e., any odd number which is less or equal $n$) theorem would give $\rk \pi _{2p-1} (G/H) = 1$. Thus, for those spaces with $n=2^{q-1}$, $q\geq 5$, Theorem~\ref{main} is not true.
\end{ex}

\subsection{The principle of the proof.} 
Since in Theorem~\ref{main}, we assume $\gg$ to be a simple compact Lie algebra, Proposition~\ref{nu=1} implies that the proof of Theorem~\ref{main}, except 
for $\gg = D_{2n}$, reduces to investigation of non-vanishing of $\de$.
According to discussion in~\ref{GT} this can be further reduced to the investigation
of non-vanishing of the differentials $\de _{\s _{i}}$ corresponding to the simple summands
in $\gg ^{'}$. This is, again by the same discussion, equivalent to the investigation on the non-decomposability
of $\rho ^{*}_{\s _{i}}$. 
 
Thus, let $\gg ^{''}$ be a simple summand in $\gg ^{'}$.
By $x_i$ $(y_i)$ we denote coordinates on the maximal Abelian subalgebra $\TT$ for $\gg$ ($\TT ^{''}$ for $\gg ^{''}$) expressing the simple roots and the Weil invariant generators for $\gg$ ($\gg ^{''}$) in classical form.
Such coordinates we call canonical coordinates below. Then in the procedure of the proof we are going to describe $\rho ^{*}_{\TT ^{''}}$ explicitely for any simple summand $\gg ^{''}$. We do it by expressing the restrictions of the canonical coordinates for $\TT$ on $\TT ^{''}$ in terms of the canonical coordinates for $\TT ^{''}$. We are able to do that because of the knowledge of the inclusions $\TT ^{''}\subset \TT$ in the case of generalised symmetric spaces. Concerning their classification we are going to consider separately generalised symmetric spaces $G/H$  with $\rk H=\rk G$  and $\rk H<\rk G$. 
   
We start with  the two sections which provide the expressions for the  canonical
coordinates of classical $\gg$ via the canonical coordinates of certain simple summand in $\gg ^{'}$ when they restrict on it. The cases of exceptional $\gg$ will be considered directly in the proof of Theorem~\ref{main}.

 \subsection{First category generalised symmetric spaces ($\rk H = \rk G$).}

From the classification of finite order automorphisms of semisimple 
Lie algebras it follows that  all the first category generalised symmetric spaces are generated by inner automorphisms of finite order~\cite{Helgason} (p.512). 
Also for any such space the semisimple part $\gg ^{'}$ of $\gg ^{\Theta}$ can be 
obtained  from the  subdiagram of the Kac-Dynkin diagram $\gg ^{(1)}$ ( so called extended Dynkin diagram). It is well known that the diagram $\gg ^{(1)}$ 
is a diagram for the extended system of simple roots 
$\{ \alpha _{0},\alpha _{1},\ldots ,\alpha _{n}\}$ for $\gg$. Thus, the
base for $\TT ^{'}$ is some subset of $\{ H_0, H_1,\ldots ,H_n \}$, where
 $H_0, H_1,\ldots ,H_n$ is a canonical base for $\TT$, i.e., the base corresponding
to the extended system of the simple roots for $\gg$~\cite{Helgason},~\cite{Kac1},~\cite{T}. These vectors are given by the conditions
 $\alpha _{i}(H_{j}) = 
\frac{2(\alpha _{i}, \alpha _{j})}{(\alpha _{j},\alpha _{j})} = a_{ij}$,
where $(a_{ij})_{0\leq i,j\leq n}$ is the matrix corresponding to the $\gg ^{(1)}$ diagram ( so called generalised Cartan matrix of a first type). 

Let $x_1,\ldots ,x_n$ be the canonical coordinates on $\TT$. There are  well known
formulas expressing the simple roots of the compact simple Lie algebras via 
the canonical coordinates~\cite{Vinberg}. Concerning that formulas we can define
matrix $P^{(1)} = (p_{ij})$, $0\leq i\leq n$, $1\leq j\leq n$ to be
\[
p_{ij} = x_{j}(H_{i}) \ . 
\]
For the simple compact Lie algebras matrix $P$ is given as follows.

1. $\gg = A_{n}$
\begin{eqnarray}\label{ap1}
p_{ii} & = & 1,\quad i = \bar{1, n} \ ,\nonumber \\
p_{i-1i} & = & -1, \quad i = \bar{1, n+1} \ , \\
p_{0 n+1} & = & 1 \ , \nonumber
\end{eqnarray}
and all the other matrix elements are equal to zero.

2. $\gg = B_{n}$
\begin{eqnarray}\label{bp1}
p_{ii} & = & 1,\quad i = \bar{1, n-1} \ , \nonumber \\
p_{i-1i} & = & -1,\quad i = \bar{1, n} \ ,  \\
p_{02} & = & -1, \quad p_{nn} = 2 \ , \nonumber
\end{eqnarray}
and all the other elements are equal to zero.

3. $\gg = C_{n}$
\begin{eqnarray}\label{cp1}
p_{ii} & = & 1,\quad i = \bar{1, n} \ , \nonumber \\
p_{i-1i} & = & -1,\quad i = \bar{1, n} \ , 
\end{eqnarray}
and all the  other elements are equal to zero.

4. $\gg = D_{n}$
\begin{eqnarray}\label{dp1}
p_{ii} & = & 1,\quad i = \bar{1, n} \ , \nonumber \\
p_{i-1i} & = & -1,\quad i = \bar{1, n} \ ,\\
p_{02} & = & -1,\quad p_{nn-1} = 1 \ ,  \nonumber
\end{eqnarray}
and all the other elements are equal to zero.

5. $\gg = G_{2}$
\begin{equation}\label{gp1}
P = \left ( \begin{array}{ccc}
-1 & 0 & 1 \\
1 & -2 & 1 \\
0 & 1 & -1   
\end{array}
\right ) 
\end{equation}

6. $\gg = F_{4}$
\begin{equation}\label{fp1}
P = \left ( \begin{array}{ccccc}
-1 & -1 & 0 & 0 \\
1 & -1 & -1 & -1 \\
0 & 0 & 0 & 2 \\
0 & 0 & 1 & -1 \\
0 & 1 & -1 & 0  
\end{array}
\right ) 
\end{equation}

7. $\gg = E_{6}$
\begin{eqnarray}\label{e6p1}
p_{ii} & = & 1,\; i = \bar{1,5}, \qquad p_{i-1i}  =  -1,\; i = \bar{2,6} \ , \nonumber \\
p_{6i} & = & \frac{-1}{2},\; i = \bar{1,3},\qquad p_{6i} = \frac{1}{2},\; i = \bar{4,6} \ , \\
p_{0\varepsilon } & = & -1, \qquad p_{6\varepsilon } = \frac{1}{2} \ ,\nonumber
\end{eqnarray}
and all the other matrix elements are equal to zero.

8. $\gg = E_{7}$
\begin{eqnarray}\label{e71}
p_{ii} & = & 1,\; i = \bar{1,6}, \qquad p_{i-1i}  =  -1,\; i = \bar{2,7} \ , \nonumber \\
p_{7i} & = & \frac{-1}{2},\; i = \bar{1,4},\qquad p_{7i} = \frac{1}{2},\; i = \bar{5,8} \ , \\
p_{07} & = & 1, \qquad p_{08} = -1 \ , \nonumber
\end{eqnarray}
and all the other element are equal to zero.

9. $\gg = E_{8}$
\begin{eqnarray}\label{e8p1}
p_{ii} & = & 1,\; i = \bar{1,7}, \qquad p_{i-1i}  =  -1,\; i = \bar{1,8} \ , \nonumber \\
x_{8i} & = & \frac{-1}{3},\; i = \bar{1,5},\qquad p_{8i} = \frac{2}{3},\; i = \bar{6,8} \ ,
\end{eqnarray}
and all the other elements are equal to zero.

Let $\gg ^{''}$ be some simple summand in $\gg ^{'}$ and $\TT ^{''}$ its maximal
Abelian subalgebra. Then the base for $\TT ^{''}$ is given by some subset of successive 
vectors in $\{ H_{0},\ldots ,H_{n} \}$.

Because of simplicity of formulation, in the following proposition we will assume
that $\TT ^{''} \neq \LL (H_{0}, H_{2}, H_{1})$ for 
 $\gg = B_{n}$  and $\TT ^{''} \neq \LL (H_{0}, H_{2}, H_{1})$, $\TT ^{''} \neq \LL (H_{n-1}, H_{n-2}, H_{n})$ for  $\gg = D_{n}$. 
\begin{prop}\label{p:canonic1}
Let $(\gg , \gg ^{\Theta}, \Theta)$ be a first category generalised symmetric algebra of a classical simple compact Lie algebra $\gg $ and $\gg ^{''}$ be some simple summand in $\gg^ {\Theta}$. Let further $x_1,\ldots ,x_n$ be the canonical coordinates on maximal Abelian subalgebra $\TT$ for $\gg$, and $y_1,\ldots ,y_l$ the canonical coordinates on maximal Abelian subalgebra $\TT ^{''}$ for $\gg ^{''}$. Then there exist $p_1, p_2,\ldots ,p_l$, $1\leq p_i\leq n$ such that the restriction of $x_1,\ldots ,x_n$  on $\TT ^{''}$ is given by
\[
x_{p_1}\to \pm y_1,\ldots ,x_{p_l}\to \pm y_l, \;\; x_{i}\to 0, \; \mbox{for} \; i\neq p_1,\ldots ,p_l \ .
\]
\end{prop}
\begin{proof}
We will proceed by considering the fixed point subalgebras for each type $\gg$.

1. $\mathbf {\gg = A_n}$, ($n\geq 2$)

The fixed point subalgebras are of the following form:
\[
\gg ^{\Theta} = \TT ^{n-\sum _{i=1}^{k}n_i}\oplus A_{n_1}\oplus \ldots \oplus A_{n_k} \ .
\]
Thus, $\gg ^{''}$ here can be only of a type $A_l$, $l\leq n$. The basis for
 $\TT ^{''}$ is
of the form 
\[
\TT ^{''} = \LL (H_{p+1},\ldots ,H_{p+l}), \;\; 0\leq p\leq n-l-1 \ ,
\]
and we may assume that it does not 
contain $H_0$, since any such  case for $\gg ^{''} = A_{l}$ is conjugate to this one by automorphism of the Kac-Dynkin diagram  $A_{n}^{(1)}$. Then, from~\eqref{ap1} follows
\[
x_{p+1}\to y_1,\ldots ,x_{p+l+1}\to y_{l+1}, \;\; x_{i}\to 0, \; i\neq p+1,\ldots ,p+l \ . 
\]

2 . $\mathbf {\gg = B_n}$, $n\geq 2$

Here the fixed point subalgebras are of the form:
\[
\gg ^{\Theta} = \TT ^{n-\sum_{i=1}^{k}n_i}\oplus D_{n_1}\oplus A_{n_2}\oplus \ldots \oplus A_{n_{k-1}}\oplus B_{n_k} \ .
\]
\begin{itemize}
\item[-] 
If $\gg ^{''} = D_l$, $l\leq n$ then
\[
\TT ^{''} = \LL (H_{l-1},\ldots ,H_{0}) \ ,
\]
and from~\eqref{bp1} and~\eqref{dp1} we get
\[ 
x_l \to -y_1,\ldots ,x_1 \to -y_l, \;\; x_{i}\to 0, \; i>l \ .
\]
\item[-] 
If $\gg ^{''} = A_l$, $l\leq n-1$, then 
\[
\TT ^{''} = \LL (H_{p+1},\ldots ,H_{p+l}), \;\; 0\leq p\leq n-l-1
\]
and from~\eqref{ap1} and~\eqref{bp1} it  follows
\[
x_{p+1}\to y_1,\ldots ,x_{p+l+1}\to y_{l+1}, \;\;  x_i \to 0, \;  i\neq p+1,\ldots ,p+l+1 \ .
\]
Note that the cases  $\TT ^{''} = \LL (H_0,H_2,\ldots ,H_{l-1})$ and $\TT ^{''} = \LL (H_1,\ldots ,H_{l-1})$, corresponding to $\gg ^{''} = A_{l}$, are the same, up to conjugation, since there is an automorphism of $B_{n}^{(1)}$ transforming one base into another.

\item[-]
If $\gg ^{''}= B_l$, $l<n$ ($l = n$ gives the trivial case) then
\[
\TT ^{''} = \LL (H_{n-l+1},\ldots , H_l) \ ,
\]
and again from~\eqref{bp1} it follows 
\[
x_{n-l+1}\to y_1, \ldots ,x_n \to y_l, \;\;  x_i \to 0, \;   i< n-l+1 \ .
\]
\end{itemize}

3. $\mathbf {\gg = C_n}$ $(n\geq 3)$

Here the fixed point subalgebras are given by:
\[
\gg ^{\Theta} = \TT ^{n-\sum _{i=1}^{k}n_i}\oplus C_{n_1}\oplus A_{n_2}\oplus \ldots \oplus  A_{n_{k-1}}\oplus C_{n_k} \ .
\]
\begin{itemize}
\item[-]
 If $\gg ^{''}=C_l$, $l\leq n$ we have two possibilities for the base of $\TT ^{''}$.

\[
 1. \;  \TT ^{''} = \LL (H_{l-1},\ldots ,H_0) \ ,
\]
 and from~\eqref{cp1} we get
\[
x_l \to -y_1,\ldots ,x_1 \to -y_l, \;\;  x_i \to 0, \; i>l \ .
\]
\[
2. \; \TT ^{''} = \LL (H_{n-l+1},\ldots ,H_n) \ ,
\]
\[
x_{n-l+1}\to y_1,\ldots ,x_n \to y_l, \;\;  x_i \to 0, \;  i< n-l+1 \ .
\]

\item[-]
If $\gg ^{''}= A_l$, $l\leq n-1$ then 
\[
\TT ^{''} = \LL (H_{p+1},\ldots ,H_{p+l}), \; 0\leq p\leq n-l-1 \ ,
\]
and from~\eqref{ap1} and~\eqref{cp1} it follows
\[
x_{p+1}\to y_1,\ldots ,x_{p+l+1}\to y_{l+1}, \;\; x_i \to 0, \; i\neq p+1,\ldots ,p+l+1 \ .
\]
\end{itemize}

4. $\mathbf {\gg = D_n}$ $(n\geq 4)$

The fixed point subalgebras are given by:

\[
\gg ^{\Theta} = \TT ^{n - \sum _{i=1}^{k}n_i}\oplus D_{n_1}\oplus A_{n_2}\oplus \ldots \oplus A_{n_{k-1}}\oplus D_{n_k} \ .
\]

\begin{itemize}
\item[-] 
If $\gg ^{''}=D_l$, $l< n$, ($l = n$ gives the trivial case) we have two possibilities for the base of $\TT ^{''}$.

\[
 1. \; \TT ^{''} = \LL (H_{l-1},\ldots ,H_{0}) \ ,
\]
then from~\eqref{dp1} it follows 
\[
x_l \to -y_1,\ldots ,x_1 \to -y_l, \;\; x_i \to 0, \;  i>l \ .
\]

\[
2. \; \TT ^{''} = \LL (H_{n-l+1},\ldots ,H_n) \ ,
\]
and again from~\eqref{dp1} it follows
\[
x_{n-l+1}\to y_1,\ldots ,x_n \to y_l, \;\;  x_i \to 0, \; i<n-l+1 \ .
\]

\item[-]
 If $\gg ^{''} = A_l$, $l\leq n-1$, then 
\[
\TT ^{''} = \LL (H_{p+1},\ldots ,H_{p+l}), \;  0\leq p\leq n-l-1 \ ,
\]
and using~\eqref{ap1} and~\eqref{dp1} we get
\[
x_{p+1}\to y_1,\ldots ,x_{p+l+1}\to y_{l+1}, \;\; x_i \to 0, \;  i\neq p+1,\ldots ,p+l+1 \ .
\]
Case $\gg ^{''} = A_{l}$ and $\TT ^{''} = \LL (H_0, H_2,\ldots ,H_{l-1})$ is conjugate by automorphism of the Kac-Dynkin $D_n^{(1)}$ to the case $\gg ^{''} = A_{l}$ and $\TT ^{''} = \LL (H_1,\ldots ,H_{l-1})$. Also the case $\TT ^{''} = \LL (H_{n-l},\ldots ,H_{n-2},H_n)$ is conjugate to the case $\TT ^{''} = \LL (H_{n-l},\ldots , H_{n-1})$.
 \end{itemize}
\end{proof}

\subsection{Second category homogeneous spaces ($\rk H < \rk G$).}

From the classification of the finite order automorphisms of the semisimple 
Lie algebras it follows that all the  second category generalised symmetric spaces are generated by an outer automorphisms of a finite order. Such  spaces exist for
$\gg = A_{2}, A_{2r}, A_{2r-1}, D_{r+1}, D_4$ and $\gg = E_6$.
In this case the semisimple part $\gg ^{'}$ of $\gg ^{\Theta}$ can be 
obtained  from the  subdiagram of the Kac-Dynkin diagrams $\gg ^{(2)}$ or $\gg ^{(3)}$.
The one-forms having $\gg ^{(2)}$ and $\gg ^{(3)}$ as the diagrams and, related to the generalised Cartan matrices of the second and third type, corresponding them vectors $\bar{H}_{i}$, $0\leq i\leq r$ on $\TT$ are given in~\cite{Helgason},~\cite{Kac1} and~\cite{T}. As in the case of the first category  generalised symmetric spaces 
we will consider matrix $P$ whose elements are defined by
$p_{ij} = x_{i}(\bar{H_j})$.

1. $\gg = A_{2r}$, $r\geq 2$ 

\begin{eqnarray}  \label{a2n}            
\bar{H_0} & = & -(H_1 +\ldots +H_{2r}) \ , \nonumber \\
\bar{H_i} & = & H_i + H_{2r-i+1}\quad (1\leq i\leq r-1) \ , \\
\bar{H_r} & = & 2(H_r + H_{r+1}) \ . \nonumber
\end{eqnarray} 

Matrix of the canonical coordinates for $A_{2r}$ in the above base  looks like:   
\begin{eqnarray}\label{a2sp2}
 p_{ii} & = & 1, \; i = \bar{1, r-1}, \qquad p_{i-1i} = -1, \; i = \bar {1, r} \ , \nonumber \\  
p_{rr} & = & 2 ,\qquad p_{jr+1} = 0, \; j = \bar {1, r} \ ,\\
 p_{ij} & = & - p_{i 2r-j+2}, \; j = \bar {1, r}, \;  i = \bar {0, r} \ , \nonumber 
\end{eqnarray}

and the other elements are equal to zero.

2. $\gg = A_{2r-1}$, $r\geq 3$

\begin{eqnarray} \label{a2n-1}
\bar{H_0} & = & -(H_1 + H_{2r-1} +2H_2 +\ldots + 2H_{2r-2}) \ ,\nonumber \\
\bar{H_i} & = & H_i + H_{2r-i}\quad (1\leq i\leq r-1) \ , \\
\bar{H_r} & = & H_r \ . \nonumber
\end{eqnarray}

Matrix of the canonical coordinates for $A_{2r-1}$ in this base is: 
\begin{eqnarray}\label{a2s-1p2}
& p_{ii} = 1, \; i = \bar{1, r},  \qquad  p_{i-1i} = -1, \; i = \bar {1, r} \ , \\
& p_{02} = -1, \qquad p_{ij} = - p_{i2r-j}, \; j = \bar{1, r} , i = \bar {0, r} \ , \nonumber
\end{eqnarray}
and the other matrix elements are equal to zero.

3. $\gg = D_{r+1}$, $r\geq 2$

\begin{eqnarray} \label{dn+1}
\bar{H_0} & = & -2(H_1 +\ldots +H_{r-1}) - (H_r + H_{r+1}) \ , \nonumber \\
\bar{H_i} & = & H_i \quad (1\leq i\leq r-1) \ ,\\
\bar{H_r} & = & H_r + H_{r+1} \ .\nonumber
\end{eqnarray}

Matrix of the canonical coordinates for $D_{r+1}$ in the above base is:
\begin{eqnarray}\label{dn+1p2}
& p_{ii} = 1, \qquad p_{ii+1} = -1, \; i = \bar{1, r-1} \ , \\ 
& p_{01} = -2, \qquad p_{rr} = 2 \ , \nonumber
\end{eqnarray}
and the other elements of the matrix are equal to zero.

4. $\gg = D_{4}$

\begin{eqnarray} \label{d4}
\bar{H_0} & = & -3H_2 - 2(H_1 + H_3 + H_4) \ ,\nonumber\\
\bar{H_1} & = & H_1 + H_3 + H_4 \ , \\
\bar{H_2} & = & H_2 \ .\nonumber
\end{eqnarray}

\begin{equation}\label{d4p2}
P = \left ( \begin{array}{cccc}
-2 & -1 & -1 & 0 \\
1 & -1 & 2 & 0 \\
0 & 1 & -1 & 0  
\end{array}
\right ) 
\end{equation}

5. $\gg = E_{6}$

\begin{eqnarray} \label{e6}
\bar{H}_{0} & = & -(2H_1 + 3H_2 + 4H_3 + 3H_4 + 2H_5 + 2H_6) \ ,\nonumber \\
\bar{H}_{1} & = & H_{1}+H_{5},\quad \bar{H}_{2} = H_{2}+H_{4} \ ,\\
\bar{H}_{3} & = & H_{3}, \quad \bar{H}_{4} =  H_{6} \ .  \nonumber
\end{eqnarray}

For proving the Theorem~\ref{main} in this case we want to  use the generators of Weil invariant polynomial algebra in the form~\eqref{takgen}. 
It is easy to see that on the vector space spanned by~(\ref{e6}) forms $a_{i}$, $b_{j}$ and $c_{ij}$ satisfy the relations
\[ 
b_{i} = - a_{7-i}, \; c_{ij} = - a_{i} + a_{7-j} \ , 
\]
and, then, the generators~\eqref{takgen} are given by
\[
I_{k} = \sum_{i=1}^{6}a_{i}^{k} + \sum_{i=2}^{5}c_{1i}^{k} + \sum _{j = 3, 4}c_{2j}^{k} \ .
\]

Thus, in order to compute the restriction of these generators on $\TT ^{'}$,
 we need to know the matrix $A$ with the elements defined by $a_{ij} = a_{j}(\bar{H}_{i})$. By straightforward computation we get that

\begin{equation}\label{e6a}
A = \left ( \begin{array}{ccccccc}
-2 & -1 & -1 & -1  & -1 & 0 \\
1 & -1 & 0 & 0 & 1 & -1 \\
0 & 1 & -1 & 1 & -1 & 0  \\
0 & 0 & 1 & -1 & 0 & 0 \\
0 & 0 & 0 & 1 & 1 & 1\\
\end{array}
\right ) \ .
\end{equation}

\begin{rem}
For $\gg = A_{2}$, we have $\gg ^{\Theta} = \TT ^{1}$ or $\gg ^{\Theta} = A_{1}$. Example~\ref{flag} and Proposition~\ref{pi34} immediately proves the
Theorem~\ref{main} in this case.
\end{rem} 

\begin{rem}\label{dodatak}
Note that~\eqref{a2n},~\eqref{a2n-1} and~\eqref{dn+1} implies that on $\TT ^{\Theta}$ the following (respectively) equalities are satisfied:
\begin{itemize}
\item $x_{r+1} = 0$, $x_{j} = - x_{2r-j+2}$, $1\leq j\leq r$ for $\gg = A_{2r}$ \ , 
\item $x_{j} = - x_{2r-j}$, $1\leq j\leq r$, for $\gg = A_{2r-1}$ \ , 
\item $x_{r+1} = 0$ for $\gg = D_{r+1}$ \ .
\end{itemize}
\end{rem}

Obviously, for the second category generalised symmetric spaces, the base for $\TT  ^{''}$ corresponding to some simple summand $\gg ^{''}$, consists of some successive vectors in $\{ \bar{H}_{0},\ldots ,\bar{H}_{r}\}$ given above.

As in the case of the first category generalised symmetric spaces, we are here also able to prove a similar statement. 
Because of the above remark, we need to consider only the canonical coordinates $x_{i}$, $1\leq i\leq r$. Again for simplicity of formulation we will assume that $\TT ^{''}\neq \LL (\bar{H}_{0}, \bar{H}_{2}, \bar{H}_{1})$ 
in a case $\gg = A_{2r}$.

\begin{prop}\label{p:canonic2}
Let $(\gg, \gg ^{\Theta}, \Theta)$ be a second category generalised symmetric algebra
of a classical compact Lie algebra $\gg$ and $\gg ^{''}$ be some simple summand in $\gg ^{\Theta}$.  Let further $x_1,\ldots ,x_r$ be the canonical coordinates on maximal Abelian subalgebra $\TT$ for $\gg$, and $y_1,\ldots ,y_l$ the canonical coordinates on maximal Abelian subalgebra $\TT ^{''}$ for $\gg ^{''}$. Then  there exist $p_1,\ldots ,p_l$, $1\leq p_i\leq r$, such that the canonical coordinates of $\gg$ restrict on $\TT ^{''}$ as follows:
\[
x_{p_1}\to \pm y_1,\ldots ,x_{p_l}\to \pm y_l, \;\;  x_{i}\to 0, \; 1\leq i\leq r, \; i\neq p_1,\ldots ,p_l \ .
\]

\end{prop}

\begin{proof}

To prove the proposition we consider separately three possible cases.

1. $\mathbf {\gg = A_{2r}}$, $r\geq 2$

Here the fixed point subalgebras are of the following form
\[
\gg ^{\Theta} = \TT ^{r-\sum _{i=1}^{k}n_i}\oplus B_{n_1}\oplus A_{n_2}\oplus \ldots \oplus
A_{n_{k-1}}\oplus C_{n_k} \ .
\]
\begin{itemize}
\item[-] 
If $\gg ^{''} = B_l$, $l\leq r$ then
\[
\TT ^{''} = \LL (\bar{H}_{r-l+1},\ldots ,\bar{H}_{r}) \ ,
\]
and from~\eqref{a2n} and ~\eqref{bp1} we get
\[
x_{r-l+1}\to y_1,\ldots ,x_r \to y_l, \;\;  x_i \to 0, \; i < r-l+1
\ .
\]

\item[-] 
If $\gg ^{''} = A_{l}$, $l\leq r-1$, then 
\[
\TT ^{''} = \LL (\bar{H}_{p+1},\ldots ,\bar{H}_{p+l}), \;  0\leq p\leq r-l-1 \ ,
\]
 and from~\eqref{a2n} and~\eqref{ap1} follows 
\[ 
x_{p+1}\to y_1,\ldots ,x_{p+l+1}\to y_{l+1}, \ ,
\]
\[
x_i \to 0, \;   i\neq p+1,\ldots ,p+l+1, \; 1\leq i\leq r \ .
\]

\item[-] 
If $\gg ^{''} = C_l$, $l\leq r$ then 
\[
\TT ^{''} = \LL (\bar{H}_{l-1},\ldots ,\bar{H}_0) \ ,
 \]
and~\eqref{a2n} and~\eqref{cp1} implies 
\[
x_l\to -y_1,\ldots ,x_1\to -y_l, \;\; x_i \to 0, \; l < i\leq r \ .
\]
\end{itemize}

2. $\mathbf {\gg = A_{2r-1}}$, $r\geq 3$

The fixed point subalgebras are given by 

\[
\gg ^{\Theta} = \TT ^{r-\sum_{i=1}^{k}n_i}\oplus D_{n_1}\oplus A_{n_2}\oplus \ldots \oplus A_{n_k}\oplus C_{n_k} \ .
\]
\begin{itemize}
\item[-] 
If $\gg ^{''} = D_l$, $l\leq r$ then 
\[
\TT ^{''} = \LL (\bar{H}_{l-1},\ldots ,\bar{H}_4,\bar{H}_3,\bar{H}_1,\bar{H}_0) \ ,
\]
and from~\eqref{a2n-1} and~\eqref{dp1} follows
\[
x_l \to -y_1,\ldots ,x_1 \to -y_l, \;\;  x_i \to 0, \;   l<i\leq r \ .
\]

\item[-] 
If $\gg ^{''} = A_l$, $l\leq r-1$
then we have 
\[
 \TT ^{''} = \LL (\bar{H}_{p+1},\ldots ,\bar{H}_{p+l}), \; 0\leq p\leq r-l-1 \ ,
\]
and again~\eqref{a2n-1} and~\eqref{ap1} implies 
\[
x_{p+1}\to y_1,\ldots ,x_{p+l+1}\to y_{l+1} \ ,
\]
\[
  x_i \to 0, \;  i\neq p+1,\ldots p+l+1, \; 1\leq i\leq r \ .
\]

We do not have to consider the case $\TT ^{''} = \LL (\bar{H}_0,\bar{H}_2,\ldots ,\bar{H}_{l-1})$, since the subalgebras $\gg ^{''} = A_{l}$ corresponding to this case and to the case $\TT ^{''} = \LL (\bar{H}_1,\bar{H}_2,\ldots ,\bar{H}_{l-1})$ are conjugate by an automorphism of the Kac-Dynkin diagram $A_{2r-1}^{(2)}$.
 
\item[-]
If $\gg ^{''} = C_{l}$, $l\leq r$ then
\[
\TT ^{''} = \LL (\bar{H}_{r-l+1},\ldots ,\bar{H}_{l}) \ , 
\]
and~\eqref{a2n-1} and~\eqref{cp1} implies
\[
x_{r-l+1}\to y_1,\ldots ,x_{r}\to y_{l}, \;\; x_{i}\to 0, \; i < r-l+1 \ .
\]
We do not have to make the difference between the cases  $\TT ^{''} = \LL (\bar{H}_{0},\bar{H}_{2},\ldots ,\bar{H}_{r})$ and  $\TT ^{''} = \LL (\bar{H}_{1},\bar{H}_{2},\ldots ,\bar{H}_{r})$ since the 
sublagebras $\gg ^{''} = C_{r}$ corresponding to these  cases are conjugated by an automorphism
of the  Kac-Dynkin diagram $A_{2r-1}^{(2)}$.

\end{itemize}

3. $\mathbf {\gg = D_{r+1}}$ 

The fixed point subalgerbra of $D_{r+1}$ has the following form
\[
\gg ^{\Theta} = \TT ^{r-\sum_{i=1}^{k}n_i}\oplus B_{n_1}\oplus A_{n_2}\oplus \ldots \oplus A_{n_{k-1}}\oplus B_{n_k} \ .
\]
\begin{itemize}
\item[-]
If $\gg ^{''} = B_l$, $l\leq r$  we have two possibilities for the base of $\TT ^{''}$.
\[
1. \; \TT ^{''} = \LL (\bar{H}_{l-1},\ldots ,\bar{H}_0) \ ,
\]
and~\eqref{dn+1} and~\eqref{bp1} implies 
\[
 x_l\to -y_1,\ldots ,x_1\to -y_l, \;\; x_i\to 0, \;  i>l \ .
\]
\[
2. \; \TT ^{''} =  \LL (\bar{H}_{r-l+1},\ldots ,\bar{H}_r) \ ,
\]
and again  we get
\[
x_{r-l+1}\to y_1,\ldots ,x_r\to y_l,\;\; x_i \to 0, \; i<r-l+1 \ .
\]

\item[-]
 If $\gg ^{''} = A_l$, $l\leq r-1$, then
\[
\TT ^{''} = \LL (\bar{H}_{p+1},\ldots ,\bar{H}_{p+l}), \;  0\leq p\leq r-l-1 \ ,
\]
and by~\eqref{dn+1} and~\eqref{ap1}
\[
x_{p+1}\to y_1,\ldots ,x_{p+l+1}\to y_{l+1}, \;\; x_i\to 0, \; i\neq p+1,\ldots ,p+l+1 \ .
\]
\end{itemize}
\end{proof}

\subsection{Proof of Theorem~\ref{main}}

In this subsection we will prove the Theorem~\ref{main} on  rational homotopy groups of generalised symmetric spaces. 

The cases not mentioned in Theorem~\ref{main} are already proved by
Proposition~\ref{thm:odd-even}.
Also, for $p = 2$ the Theorem is proved by Proposition~\ref{pi34}. So, in the proof we consider the common exponents different from $2$.

\begin{proof}
To prove  the Theorem we will consider separately when $\gg$ is of classical and exceptional type  and we will also differ the cases of the first and the second category generalised symmetric spaces.

1. {\em $G$ is of a classical type and $\rk G=\rk H$}\\
Let $\gg ^{''}$ be some simple summand in $\gg ^{\Theta}$ and $\TT ^{''}$ its maximal 
Abelian subalgebra ($\TT ^{''}\neq \LL (H_{0}, H_{2}, H_{1})$ for $\gg = B_{n}$ or $\gg = D_{n}$ and $\TT ^{''}\neq \LL (H_{n-1}, H_{n-2}, H_{n})$ for $\gg = D_{n}$). Let $x_1,\ldots ,x_n$ be the canonical coordinates on $\TT$ and
$y_1,\ldots ,y_l$ the canonical coordinates on $\TT ^{''}$. 
Using  Proposition~\ref{p:canonic1}, we get that the restriction of $x_1,\ldots x_n$ on $\TT ^{''}$  has the following form:
\[
x_{p_i}\to \pm y_i, \;\; 1\leq i\leq l \quad x_{i}\to 0, \;\; i\neq p_1,\ldots ,p_l 
\ .
\]
 
Therefore, the restrictions on $\TT ^{''}$ of the Weil invariant generators
 for the classical groups given by~\eqref{agen} -~\eqref{takgen} have the form:
\[
 \rho ^{*}_{\TT ^{''}}(P_{k}) = (\pm y_1)^k+\ldots +(\pm y_l)^k \ ,
\]
\[
\rho ^{*}_{\TT ^{''}}(P^{'}_{n}) = y_{1}\cdots ,y_{n} \quad \mbox{for} \quad \gg ^{'} = \gg ^{''} = A_{n-1} \ ,
\]
\[
\rho ^{*}_{\TT ^{''}}(P^{'}_{n}) = 0 \quad \mbox{for} \quad \gg ^{''} \neq A_{n-1} \ .
\]  
Thus, if $p = k_i = l_j$ is the common exponent for  $\gg $ and $\gg ^{''}$, we get 
\begin{equation}\label{P_k}
\de _{\TT ^{''}}(z_{p}) = (\pm y_1)^p+\ldots +(\pm y_l)^p\neq 0 \ .
\end{equation}
For $p = n$ and $\gg = D_{n}$ we have
\begin{equation}\label{P'A}
\de (z^{'}_{n}) = y_{1}\cdots y_{n}\neq 0 \quad \mbox{for} \quad \gg ^{'} = \gg ^{''} = A_{n-1} \ ,
\end{equation}
\begin{equation}\label{P'}
\de (z_{n}^{'}) = 0 \quad \mbox{for} \quad \gg ^{''}\neq A_{n-1}
\end{equation}

Note that, when $\gg = D_{n}$ and $n$ is odd, $n$ can be a common exponent only when
$\gg ^{'} = \gg ^{''} = A_{n-1}$.

For $\gg ^{''} = A_3$ and $\TT ^{''} = \LL (H_{0}, H_{2}, H_{3})$, which we  have for $\gg = B_n$ and $\gg = D_{n}$, and more $\TT ^{''} = \LL (H_{n-2}, H_{n-1}, H_{n})$ for  $\gg = D_n$,
the only common exponent for $\gg ^{''}$ and $\gg$ is 4.
 From~\eqref{a2n-1} and~\eqref{bp1}, and,~\eqref{a2n-1} and~\eqref{dp1} we get  
\[
\rho ^{*}_{\TT ^{''}}(P_{4}) = (y_1+y_2)^{4}+(y_1+y_3)^{4}+(y_2+y_3)^{4} \ ,
\]
and it is obvious that $\rho ^{*}_{\TT ^{''}}(P_{4})$ can not be written in the form $c\dot (y_1^2+y_2^2+y_3^2+y_{1}y_{2}+y_{1}y_{3}+y_{2}y_{3})^2$. It follows that $\de (z_{4})\neq 0$.

$\bullet$ Thus, if $(\gg, p)\neq (D_{2n}, 2n)$, we have $\nu (k_i) = 1$ and from ~\eqref{P_k},~\eqref{P'A} and Proposition~\ref{nu=1} it follows
\[
\dim \pi _{2p}(G/H)\otimes \Q = \nu (l_{j}) -1, \; \pi _{2p-1}(G/H)\otimes \Q = 0 \ .
\]
$\bullet$ For $(\gg ,p) = (D_{2n}, 2n)$ we have $\nu (k_i) = 2$, and, for $\gg ^{'}\neq A_{n-1}$ from~\eqref{P_k} and~\eqref{P'} it follows
\[
\dim \pi _{4n}(G/H)\otimes \Q = \nu (l_{j}) - 1, \qquad \pi _{4n-1}(G/H)\otimes \Q = \Q \ .
\]
$\bullet$ For $\gg = D_{2n}$ and $\gg ^{'} = A_{n-1}$ using~\eqref{P_k} and~\eqref{P'A} we get
\[
\dim \pi _{4n}(G/H)\otimes \Q = \nu (l_{j}) - 1, \qquad \pi _{4n-1}(G/H) \otimes \Q = 0 \ .
\] 

2. {\em $G$ is of a classical type and $\rk H<\rk G$}\\
$\bullet$ Let $p = 2k-1$. 
For $\gg = A_{2r}$ or $\gg = A_{2r-1}$ from Remark~\ref{dodatak}
it follows that $\rho ^{*}(P_{2k-1}) = 0$, what implies $\de (z_{2k-1}) =0$. 
For $\gg = D_{r+1}$, the only odd exponent can be $r+1$ in the case  when $r$ is even. 
Concerning formula for $\gg ^{\Theta}$ in this case, given in the proof of the Proposition~\ref{p:canonic2}, we see that $r+1$ can not be exponent for $\gg ^{\Theta}$. Again using Proposition~\ref{nu=1}  we get in this case 
\[
\dim \pi _{2p}(G/H)\otimes \Q = \nu (l_j), \;\;  \pi _{2p-1}(G/H)\otimes \Q = \Q \ .
\]
   
$\bullet$ Let  $ p = 2k$ be a common exponent for $\gg$ and $\gg ^{''}$.
From Proposition~\ref{p:canonic2}  we get that $\rho ^{*}_{\TT ^{''}}$  ($\TT ^{''}\neq \LL (\bar{H}_{0}, \bar{H}_{2}, \bar{H}_{1})$ for  $\gg = A_{2r}$) for $\gg$ of a type
$A_{2r}$ or $A_{2r-1}$ is given by
\[
\rho ^{*}_{\TT ^{''}}(P_{2k}) = 2(y_1^{2k} + \ldots + y_l^{2k}) \ ,
\]
and therefore $\de (z_{2k})$ is non-trivial. For $\gg^ {''} = A_{3}$ and $\TT ^{''} = \LL (\bar{H}_{0}, \bar{H}_{2}, \bar{H}_{1})$  which we may have in  the case $\gg = A_{2r}$, from~\eqref{a2sp2} we get
\[
\rho ^{*}_{\TT ^{''}}(P_{4}) = (y_{1}+y_{2})^{4} + (y_{1}+y_{3})^{4} + (y_{2}+y_{3})^{4}
\]
and as in a previous such a case it follows that $\de (z_{4})\neq 0$.

For $\gg = D_{r+1}$ and $2k\neq r+1$, using Proposition~\ref{p:canonic2} it follows
\[
\rho ^{*}_{\TT ^{''}}(P_{2k}) = y_1^{2k} + \ldots + y_l^{2k} \ , \ 
\]
and therefore $\de (z_{2k})\neq 0$.
Since in all these cases $\nu (k_{i}) = 1$, Proposition~\ref{nu=1} implies
\[
\dim \pi _{2p}(G/H)\otimes \Q = \nu (l_j) -1 
\qquad \pi _{2p -1}(G/H)\otimes \Q = 0 \ .
\]

For $\gg = D_{r+1}$ and $r+1 = 2k$ we have $\nu (k_i) = 2$ and Proposition~\ref{p:canonic2} and Remark~\ref{dodatak} give $\de (z_{2k})\neq 0$
and $\de (z^{'}_{2k}) = 0$. So, we get
\[
\dim \pi _{2p}(G/H)\otimes \Q = \nu (l_{j}) -1, \qquad \pi _{2p-1}(G/H)\otimes \Q = \Q \ .
\]

3. $\gg = D_4$ {\it and} $\rk H = 2$

Any automorphism of $D_{4}$ giving this case one can assume to be generated by 
the third order automorphism of the Dynkin diagram for $D_{4}$, see~\cite{Helgason},~\cite{Kac1}. 

Let $x_1, x_2, x_3, x_4$ be the canonical coordinates on $\TT$ for $D_4$. Then on any $\TT ^{'}$ belonging to this case,~\eqref{d4p2} implies $x_4 = 0$ and $x_1 = x_2+x_3$.
Thus $\rho ^{*}_{t^{'}}(P_{4}^{'}) =0$ and $\rho ^{*}_{t^{'}}(P_{4}) = \frac{1}{4}(\rho ^{*}_{t^{'}}(P_{2}))^{2}$. Only interesting simple summand  is $G_2$ having $6$ as a common exponent $(\neq 2)$ with $D_{4}$. Since, for it, $\rho ^{*}(P_{2})$ and $\rho ^{*}(P_{6})$ are functionally independent ( because of finites of cohomology), it follows 
that $\rho ^{*}(P_{6})$ has to contain a Weil invariant generator $Q_6$ of $G_2$. Thus, $\de (z_{6})\neq 0$ and $\pi _{12}(D_{4}/G_{2})\otimes \Q =\pi _{11}(D_{4}/G_{2})\otimes \Q = 0$.

4. $\gg = G_2$

There is no generalised symmetric spaces of $G_2$ whose fixed point subgroup has the exponent $6$.

5. $ g = F_4$

All exponents of $F_{4}$ have multiplicity $1$. The simple summands in $\gg ^{'}$ having the common exponents with $\gg $  are $B_3$, $B_4$ and $C_4$. 
\begin{itemize}
\item
 If $\gg ^{''} = B_3$, then $k_{2} = l_{3} = 6$ and
\[
\TT ^{''} = \LL (H_4, H_3, H_2) \ .
\]
From~\eqref{bp1} and~\eqref{fp1} we get
\[
x_1\to 0, \; x_2\to y_1,\; x_3\to y_2,\; x_4\to y_3 \ .
\]
Using the above formulas we can calculate $\rho ^{*}_{\TT ^{''}}(P_{6})$, and, computing it at the
points $(1,0,0)$, $(1,1,0)$ and $(1,1,1)$, we conclude that it can not be written in the form $c_{1}Q_{2}^{3} + c_{2}Q_{2}Q_{4}$. Thus, $\de (z_{6})\neq 0$ and Proposition~\ref{nu=1} implies the Theorem.

\item
If  $\gg ^{''} = B_4$, then $k_{2} = l_{3} = 6$, $k_{3} = l_{4} = 8$ and 
\[
\TT ^{''} = \LL (H_0, H_4, H_3, H_2) \ .
\]
From~\eqref{bp1} and~\eqref{fp1} it follows
\[
x_1\to -y_1, \;  x_2\to y_2, \;  x_3\to y_3, \; x_4\to y_4 \ .
\]
We prove that $\de (z_{6})\neq 0$ arguing as in the previous case \\(the forth coordinate will be zero). To prove that $\rho ^{*}_{\TT ^{''}}(P_{8})$ 
contains $Q_8$, we can just compute $\rho ^{*}_{\TT ^{''}}(P_{8})$ at the points $(1,0,0,0)$, $(1,1,0,0)$, $(1,1,1,0)$, $(1,1,1,1)$ to conclude that
it can not be written in the form $c_{1}Q_{2}^{4}+c_{2}Q_{2}^{2}Q_{4}+c_{3}Q_{2}Q_{6}+c_{4}Q_{4}^{2}$. Then from Proposition~\ref{nu=1} the Theorem follows. 

\item
If  $\gg ^{''} = C_3$, then $k_{2} = l_{3} = 6$ and 
\[
\TT ^{''} = \LL (H_1, H_2, H_3) \ .
\]
Using~\eqref{fp1} and~\eqref{cp1} we get
\[
x_1\to y_1, \; x_2\to -y_1, \; x_3\to y_2+y_3, \; x_4\to y_2-y_3 \ .
\]
Computing $\rho ^{*}_{\TT ^{''}}(P_{6})$ at the points $(1,0,0)$, $(1,1,0)$ and $(1,1,1)$ we get that $\de (z_{6})\neq 0$.
\end{itemize}

6. $\gg = E_6$

All exponents for $E_{6}$ have multiplicity $1$.

a) $\rk G = \rk H$

The simple summands in $\gg ^{'}$ having the common exponents with  $E_6$  are $A_4$, $A_5$, $D_{4}$ and $D_5$.

\begin{itemize}
\item 
If $\gg ^{''} = A_4$, then $k_{2} = l_{4} = 5$,
\[
\TT ^{''} = \LL (H_1, H_2, H_3, H_4) \ ,
\]
 and by~\eqref{e6p1} and~\eqref{ap1} we get
\[
x_i\to y_i,\; 1\leq i\leq 5, \; \;  x_6,\varepsilon \to 0 \ .
\]
Any other case for $\gg ^{''} = A_{4}$ is conjugate to this one by an automorphism of the Kac-Dynkin diagram $E_6^{(1)}$.

Using the above formulas we get
\[
\rho ^{*}_{\TT ^{''}}(P_5) = 2\sum _{i=1}^{5}y_i^{5} - \sum _{1\leq i<j\leq 5}^{4}(y_i+y_j)^{5} \ ,
\]
We need to prove that the expression for $\rho ^{*}_{\TT ^{''}}(P_{5})$ via the generators of the Weil
invariants for $A_4$ contains $Q_5$. Let us consider polynomial $Q = c\cdot
Q_{2}Q_{3}$. Then $\rho ^{*}_{\TT ^{''}}(P_5)(1,0,0,0,0) = -2$,  $\rho ^{*}_{\TT ^{''}}(P_5)(1,1,0,0,0) = -34$, while $Q(1,0,0,0,0) = c$,  $Q(1,0,0,0,0) = 4c$, what implies that $\rho ^{*}_{\TT ^{''}}(P_5)$ can not be written  in the form   $c\cdot Q_{2}Q_{3}$. Thus, $\de (z_{5})\neq 0$, what by Proposition~\ref{nu=1} proves the  the Theorem in this case.
\item 
If $\gg ^{''} = A_5$ then $k_{2} = l_{4} = 5$ and $k_{3} = l_{5} = 6$ and,
\[
\TT ^{''} = \LL (H_1, H_2, H_3, H_4, H_5) \ ,
\]
 up to conjugation.

Then~\eqref{ap1} and~\eqref{e6p1} implies
\[
x_i\to y_i, \; 1\leq i\leq 6, \;\;  \varepsilon \to 0 \ .
\]
 To prove that $\de (z_{5})\neq 0$  we argue as in the previous case. Now, $\rho ^{*}_{\TT ^{''}}(P_6) = 2\sum _{i=1}^{6}y_i^{6} + \sum_{1\leq i<j\leq 6}(y_{i}+y_{j})^{6}$, and computing its values at the points $(1,0,0,0,0,0)$, $(1,1,0,0,0,0)$ and $(1,1,1,0,0,0)$, we conclude that it can not be written in the form $c_{1}Q_{2}^{3}+c_{2}Q_{2}Q_{4}+c_{3}Q_{3}^{2}$. Thus, $\de (z_{6})\neq 0$.

\item
If $\gg ^{''} = D_4$, then $k_{3} = l_{3} = 6$,
\[
\TT ^{''} = \LL (H_2,H_3,H_4,H_6) \ ,
\]
and~\eqref{e6p1} and~\eqref{dp1} implies 
\[
x_1\to -y, \;\; x_2\to y_1 - y, \;\; x_3\to y_2 - y,
\]
\[ 
x_4\to y_3 - y,\;\; x_5\to y_4 - y,\;\; x_6\to y, \;\;  
\varepsilon\to y \ ,
\]
where $y = \frac{1}{4}\sum_{i=1}^{4}y_i$.
Using the above formulas, we can calculate $\rho ^{*}_{\TT ^{''}}(P_{6})$ at the points $(1,0,0,0)$, $(1,1,0,0)$, $(1,1,1,1)$ and conclude that it can not be written in the form $c_{1}Q_{2}^{3}+c_{2}Q_{2}Q_{4}$, i.e., $\de (z_{6})\neq 0$.

 \item
 If $\gg ^{''} = D_5$, then $k_{2} = l_{5} = 5$, $k_{3} = l_{3} = 6$, $k_{4} = l_{4} = 8$, 
\[
\TT ^{''} = \LL (H_1, H_2, H_3, H_4, H_6)  \ ,
\]
and we have 
\[
x_1\to y_1 - y, \;\; x_2\to y_2 - y, \;\;  x_3\to y_3 - y \ ,
\]
\[
x_4\to y_4 - y , \;\; x_5\to y_5 - y , \;\; 
x_6\to y, \;\;
\varepsilon \to y \ ,
\]
where $y = \frac{1}{4}\sum_{i=1}^{5}y_i$.
Any other case for $\gg ^{''} = D_{5}$ is conjugate to this one by an automorphism of the Kac-Dynkin diagram $E_6^{(1)}$. Since $l_5 = 5$ is the only odd exponent for $D_{5}$, obviously, $\de (z_{5}) = 0$ if and only if 
$\rho ^{*}_{\TT ^{''}}(P_{5}) = 0$. Using the above 
formulas it is easy to conclude that $\rho ^{*}_{\TT ^{''}}(P_{5})\neq 0$.
We prove that $\de (z_{6})\neq 0$ as in the previous case. For $P_8$, calculating its restriction at the points $(1,0,0,0,0)$, $(1,1,0,0,0)$, $(1,1,1,0,0)$ and $(1,1,1,1,0)$, we get that $\rho ^{*}_{\TT ^{''}}(P_{8})$ can not be written in the form $c_{1}Q_{2}^{4}+c_{2}Q_{2}^{2}Q_{4}+c_{3}Q_{2}Q_{6}+c_{4}Q_{4}^{2}$,  and thus $\de (z_{8})\neq 0$.
\end{itemize}

b) $\rk H<\rk G$

The simple summands in $\gg ^{'}$ having the common exponents 
with $E_6$  are $F_4$, $B_3$, $C_3$ and $C_4$. We will use here the generators of the Weil invariant polynomial algebra for $E_6$ in the form~\eqref{takgen}.
By straightforward calculation we get the following.
\begin{itemize}
\item[-]
 If $\gg ^{''} = B_3$, then $k_{3} = l_{3} = 6$,
\[
\TT ^{''} = \LL  (\bar{H}_{4}, \bar{H}_{3}, \bar{H}_{2}) \ ,
\]
and using~\eqref{e6a} and~\eqref{bp1} we get
\[
a_{1}\to 0,\;\; a_{2}\to \frac{1}{2}(y_1+y_2+y_3), \;\;  a_{3}\to \frac{1}{2}(y_1+y_2-y_3) \ ,
\]
\[
a_{4}\to \frac{1}{2}(y_1-y_2+y_3), \;\; a_{5}\to \frac{1}{2}(y_1-y_2-y_3),  \;\; a_{6}\to y_1 \ .
\]

If we calculate the restriction of $I_{6}$ on $\TT^{''}$ at the points $(1,0,0)$, $(1,1,0)$, $(1,1,1)$ we can easily see that it can not be written in the form $c_{1}Q_{1}^{3}+c_{2}Q_{1}Q_{2}$. Thus, $\de (z_{6})\neq 0$.
 
\item[-] 
If $\gg ^{''} = C_3$, then $k_{3} = l_{3} = 6$,
\[
\TT ^{''} = \LL (\bar{H}_{1}, \bar{H}_{2}, \bar{H}_{3}) \ ,
\]
 and using~\eqref{e6a} and~\eqref{cp1} we get :
\[
a_{1}\to y_1, \;\; a_2\to y_2, \;\; a_3\to y_3 \ ,
\]
\[a_4\to -y_3, \;\; a_5\to -y_2, \;\; a_6\to -y_1 \ . 
\]

Computing $\rho ^{*}_{\TT ^{''}}(I_{6})$ at the points $(1,0,0)$, $(1,1,0)$ and $(1,1,1)$ we conclude that it can not be written in the form $c_{1}Q_{1}^{3}+c_{2}Q_{1}Q_{2}$.

\item[-]
If $\gg ^{''} = C_4$, then $k_{3} = l_{3} = 6$ and $k_{4} = l_{4} = 8$,
\[
\TT ^{''} = \LL (\bar{H}_{0}, \bar{H}_{1}, \bar{H}_{2}, \bar{H}_{3}) \ ,
\]
and for the forms given in~\ref{takgen} we get
\[
a_1\to -y_1+y_2, \;\; a_2\to -y_1+y_3, \;\; a_3\to -y_1+y_4 \ ,
\]
\[
a_4\to -y_1-y_4, \;\; a_5\to -y_1-y_3, \;\; a_6\to -y_1-y_2 \ .
\]

As in the previous cases we can easily show that $\rho ^{*}_{\TT ^{''}}(I_{6})$, $\rho ^{*}_{\TT ^{''}}(I_{8})$   contain
$Q_6$ and $Q_8$ (respectively) as the linear parts.

\item[-]
For $\gg ^{'} = F_{4}$, in~\cite{Takeuchi} is proved that $H^{*}(E_{6}/F_{4}) = \wedge (z_{5}, z_{7})$. Thus, the Cartan-Serre theorem proves that the Theorem is true.
\end{itemize}

7. All generalised symmetric algebras with  $\gg = E_7$ and $\gg = E_8$ are of the first category. 
Concerning the exponents for $E_7$ and $E_8$  the interesting cases for $\gg ^{''}$ we need to consider  are $A_7$, $A_6$, $A_5$, $D_6$, $D_5$, $D_4$, $E_6$ for $E_7$ and $A_8$, $A_7$, $D_8$, $D_7$, $D_6$, $D_5$, $E_7$, $E_{6}$ for $E_8$. Using the generators of the Weil invariant polynomial algebras for $E_7$ and $E_8$ given in~\cite{MT}, the proof goes in the same way as in the previous cases.
Being technical in its nature we omit it here.
\end{proof}

In particular, Theorem~\ref{main}, which we just proved, gives us the rational homotopy 
groups of an irreducible simply connected compact symmetric spaces. We list them in Appendix~\ref{A} that follows.
\begin{rem} 
For the spaces denoted in the table by $(^{\dagger})$ we assume that $k\leq \frac{n+1}{2}$. For those denoted by $(^{\ddagger})$ we assume that $k\leq \frac{n}{2}$, and by $(^{\diamond})$ that $k\leq \frac{n-1}{2}$. 
\end{rem}
\begin{rem}
Note that, except in the cases denoted in the table by $(^{\star})$ for $k$ even, by $(^{\circ})$ for $n-k$ even and 
by both for $k = \frac{n}{2}$, the  dimensions of all non-trivial rational homotopy groups are $1$. The dimension is $2$ in the following cases.
\begin{itemize}
\item[-] Denoted by$(^{\star})$ for $k$ even and $k\neq \frac{n}{2}$; 
\item[-] Denoted by $(^{\circ})$ for $n-k$ even and $k\neq \frac{n}{2}$;
\item[-] For $k$ odd and  $k = \frac{n}{2}$.
\end{itemize}
For $k$ even and $k = \frac{n}{2}$, it is $3$.
\end{rem}
\newpage

\small{
\begin{appendix}
\section{Rational homotopy groups of irreducible simply connected compact symmetric spaces}\label{A}

\begin{tabular}{|l||l|} \hline
\hspace{2cm}M & \hspace{1cm} Non-trivial $\pi _{q} \otimes \Q$ \\ \hline \hline
$SU(2n+1)/SO(2n+1)$ & 
$\begin{array}{ll} q = 2p-1 & p = 3, 5,\ldots ,2n+1 
\end{array}$ 
\\ \hline
$SU(2n)/SO(2n)$ &  
$\begin{array}{ll} q = 2n &  \\
q = 2p-1 & p = 3,\ldots, 2n-1, 2n 
\end{array}$  
\\ \hline
$SU(2n)/Sp(n)$ & 
$\begin{array}{ll} q = 2p-1 & p = 3, 5,\ldots, 2n-1 
\end{array}$
\\ \hline
$^{\dagger}U(n+1)/U(k)\times U(n-k+1)$ & 
$\begin{array}{ll} q = 2p &  p = 1,2,\ldots ,k \\ q = 2p-1 & p = n-k+2,\ldots ,n+1 
\end{array}$ 
\\ \hline
$SO(2n+1)/SO(2)\times SO(2n-1)$ & 
$\begin{array}{l}
q = 2, 4n-1 \end{array}$
\\ \hline
$^{\dagger}SO(2n+1)/SO(2k)\times SO(2n+1-2k)$ & 
$\begin{array}{lll} 
q = 2k^{\star} & & \\
q = 4p &  p = 1,\ldots ,k-1 & \\
q = 4p-1 & p = n-k+1,\ldots , n & 
\end{array}$ 
\\ \hline
$SO(2n+1)/SO(2n) = S^{2n}$ & 
$\begin{array}{ll} q = 2n & \\ q = 4n-1 & \
\end{array}$ 
\\ \hline
$Sp(n)/U(n)$ & 
$\begin{array}{ll} q = 2p & p = 1,3,\ldots ,2[\frac{n-1}{2}]+1  \\
q = 4p-1 &  p =  [\frac{n}{2}]+1\ldots ,n  
\end{array}$ 
\\ \hline
$^{\ddagger}Sp(n)/Sp(k)\times Sp(n-k)$ & 
$\begin{array}{ll} q = 4p & p = 1,\ldots ,k \\
q = 4p-1 &  p = n-k+1,\ldots ,n  
\end{array}$ 
\\ \hline
$SO(2n)/SO(2)\times SO(2n-2)$ & 
$\begin{array}{l}
q = 2, 2n-2, 4n-5, 2n-1 
\end{array}$ 
\\ \hline
$^{\ddagger}SO(2n)/SO(2k)\times SO(2n-2k)$ &
$\begin{array}{lll} q = 2k^{\star}, 2n-2k^{\circ}, 2n-1 & &\\
q= 4p & p = 1,2,\ldots, k-1 & \\
q = 4p-1 & p = n-k,\ldots ,n-1 & 
\end{array}$
\\ \hline
$^{\diamond}SO(2n)/SO(2k+1)\times SO(2n-2k-1)$ &
$\begin{array}{lll}
q = 2n-1 & & \\
q = 4p & p = 1,\ldots ,k & \\
q = 4p-1 &  p = n-k,\ldots ,n-1 &  
\end{array}$
\\ \hline
$SO(2n)/SO(2n-1) = S^{2n-1}$ &
$\begin{array}{l}
q = 2n-1
\end{array}$
\\ \hline
$SO(2n)/U(n)$ & 
$\begin{array}{ll} 
q = 2p & p = 1, 3, 5,\ldots ,2[\frac{n-1}{2}]+1 \\
q = 4p-1 & p = [\frac{n}{2}]+1,\ldots ,n-1  
\end{array}$
\\ \hline
$G_{2}/SO(4)$ &
$\begin{array}{l}
q = 8, 11
\end{array}$
\\ \hline
$F_{4}/SU(2)\cdot Sp(3)$ &
$\begin{array}{l}
q = 4, 8, 15, 23
\end{array}$
\\ \hline
$F_{4}/Spin(9)$ &
$\begin{array}{l}
q = 8, 23
\end{array}$
\\ \hline
$E_{6}/PSp(4)$ & 
$\begin{array}{l} 
q = 8, 9, 17, 23
\end{array}$
\\ \hline
$E_{6}/F_{4}$ &
$\begin{array}{l}
q = 9, 17 
\end{array}$
\\ \hline
$E_{6}/SU(2)\cdot SU(6)$ & 
$\begin{array}{l}
q = 4, 6, 8, 14, 15, 17, 23
\end{array}$
\\ \hline
$AdE_{6}/T^{1}\cdot Spin(10)$ & 
$\begin{array}{l}
q = 2, 8, 17, 23
\end{array}$
\\ \hline
$E_{7}/SU^{*}(8)$ &
$\begin{array}{l}
q = 6, 8, 10, 14, 19, 23, 27, 35
\end{array}$
\\ \hline
$E_{7}/SU(2)\cdot Spin(12)$ & 
$\begin{array}{l}
q = 4, 8, 12, 23, 27, 35 
\end{array}$ 
\\ \hline
$AdE_{7}/T^{1}\cdot E_{6}$ & 
$\begin{array}{l}
q = 2, 10, 18, 19, 27, 35
\end{array}$ 
\\ \hline
$E_{8}/SO(16)$ &
$\begin{array}{l} 
q = 8, 12, 16, 20, 35, 39, 47, 59
\end{array}$
\\ \hline
$E_{8}/SU(2)\cdot E_{7}$ & 
$\begin{array}{l}
q = 4, 12, 20, 39, 47, 59 
\end{array}$
\\ \hline
\end{tabular}

\section{Generators of the Weil invariant polynomial algebras for the simple compact Lie algebras}\label{B}

\begin{itemize}
\item $\gg = A_{n}$, $(n\geq 1)$

\begin{equation} \label{agen}
 P_{k} = \sum_{i=1}^{n+1}x_{i}^{k},\quad  k = 2, 3, 4,\ldots ,n+1\quad (\sum_{i=1}^{n+1}x_{i} = 0)  \ ;
\end{equation}

\item $\gg = B_{n}$, $(n\geq 2)$, $\gg = C_{n}$, $(n\geq 3)$

\begin{equation}\label{cgen}
P_{k} = \sum_{i=1}^{n}x_{i}^{k},\quad k = 2, 4, 6,\ldots ,2n \ ;
\end{equation}
 
\item $\gg = D_{n}$, $(n\geq 4)$

\begin{equation}\label{dgen}
P_{k} = \sum_{i=1}^{n}x_{i}^{k},\quad k = 2, 4,\ldots , 2n-2 \quad P_{n}^{'} = x_{1}\cdots x_{n} \ ;
\end{equation}

\item $\gg = G_{2}$
\begin{equation}\label{g2gen}
P_{2} = 2\sum_{j=1}^{3}x_{j}^{2},\quad P_{6} = 2\sum_{j=1}^{3}x_{j}^{6}\quad  (\sum_{i=1}^{3}x_{i} = 0) \ ;
\end{equation}

\item $\gg = F_{4}$
\begin{equation}\label{f4gen}
P_{k} = \sum_{i=1}^{4}(\pm x_{i})^{k} + \{ \frac{1}{2}(\pm x_1\pm x_2\pm x_3\pm x_4)\}^{k},\quad k = 2, 6, 8, 12 \;
\end{equation}

 \item $\gg = E_{6}$

Standard generators~\cite{cox}:

\begin{equation}\label{e6gen}
P_{k} = \sum _{i=1}^{6}(x_{i}\pm \varepsilon )^{k} + \sum _{1\leq i<j\leq 6}(-x_{i}-x_{j})^{k} \;\; k = 2, 5, 6, 8, 9, 12 \ .
\end{equation}

Generators given in~\cite{Takeuchi}:
\begin{equation}\label{takgen}
I_{k} =  \frac{1}{2}(\sum_{i=1}^{6}a_{i}^{k} + \sum_{i=1}^{6}b_{i}^{k} + 
\sum_{i, j=1}^{6}c_{ij}^{k})\quad k = 2, 5, 6, 8, 9, 12 \ ,
\end{equation}
\
 where $a_i$, $b_i$ and $c_{ij}$ are given by
\[
a_{i} = x_{i}+\frac{1}{2}\sum_{i=1}^{6}x_{i}, \; \; 
b_{i} = x_{i} - \frac{2}{3}\sum_{i=1}^{6}x_{i}, \; 1\leq i\leq 6 \ ,
\]
\[
c_{ij} = -x_{i}-x_{j}+\frac{1}{2}\sum_{i=1}^{6}x_{i}, \; 1\leq i<j\leq 6 \ .
\]

\end{itemize}

Expressions for the generators of the Weil invariant polynomial algebras for $E_7$ and $E_8$ are technically complicated and we omit them here. 
\end{appendix}
}

\bibliographystyle{amsplain}

\bigskip
\end{document}